\documentclass{interact}

\usepackage[utf8]{inputenc}
\usepackage{lmodern}
\usepackage[T1]{fontenc}
\usepackage{microtype}
\usepackage{csquotes}

\usepackage{enumitem}

\usepackage{amsmath}
\usepackage{amssymb}
\usepackage{amsthm}
\usepackage{mathtools}
\usepackage{dsfont}
\usepackage{bm}
\usepackage{multirow}

\usepackage[svgnames]{xcolor}
\usepackage{tikz}
\usetikzlibrary{positioning,calc,matrix}

\DeclarePairedDelimiter{\norm}{\lVert}{\rVert}

\DeclareMathOperator{\image}{im}
\DeclareMathOperator{\rank}{rk}

\DeclareMathOperator{\diag}{diag}

\usepackage{hyperref}
\usepackage{apacite}
\usepackage[nameinlink]{cleveref}

\newtheorem{theorem}{Theorem}\Crefname{theorem}{Theorem}{Theorems}
\newtheorem{proposition}[theorem]{Proposition}\Crefname{proposition}{Proposition}{Propositions}
\Crefname{corollary}{Corollary}{Corollaries}
\newtheorem{lemma}[theorem]{Lemma}\Crefname{lemma}{Lemma}{Lemmas}

\theoremstyle{definition}
\newtheorem{definition}[theorem]{Definition}\Crefname{definition}{Def.}{Defs.}
\newtheorem{remark}[theorem]{Remark}\Crefname{remark}{Remark}{Remarks}
\Crefname{conjecture}{Conjecture}{Conjectures}
\newtheorem{example}[theorem]{Example}\Crefname{example}{Example}{Examples}
\Crefname{property}{Property}{Properties}
\newtheorem{assumption}[theorem]{Assumption}\Crefname{assumption}{Assumption}{Assumptions}
\newtheorem*{algorithm*}{Algorithm}\Crefname{algorithm}{Algorithm}{Algorithms}

\newcommand{\cIV}{\mathfrak{C}_{[E,A,B]}}
\newcommand{\setdef}[2]{\left\{\ #1\ \left|\ \vphantom{#1} #2 \right.\right\}}

\begin{document}

\title{Model predictive control for singular differential-algebraic equations}
\author{Achim Ilchmann\(^1\), Jonas Witschel\(^1\), Karl Worthmann\(^1\)}
\thanks{\(^1\)Contact: \{achim.ilchmann,jonas.witschel,karl.worthmann\}@tu-ilmenau.de, TU Ilmenau, Ehrenbergstraße 29, 98693 Ilmenau, Germany}

\maketitle

\begin{abstract}
We study model predictive control for singular differential-algebraic equations with higher index. This is a novelty when compared to the literature where only regular differential-algebraic equations with additional assumptions on the index and/or controllability are considered. By regularization techniques, we are able to derive an equivalent optimal control problem for an ordinary differential equation to which well-known model predictive control techniques can be applied. This allows the construction of terminal constraints and costs such that the origin is asymptotically stable w.r.t.\ the resulting closed-loop system.
\end{abstract}

\section{Introduction}

Differential-algebraic equations~(DAEs) play an important role in the modelling of electrical networks, mechanical multi-body systems or chemical plants~\cite{campbellapplications}. To derive the system model, often automatic modelling techniques are employed~\cite{Riaz08}. This may lead to singular differential-algebraic systems with higher index which pose special challenges w.r.t.\ control.

We want to control DAEs using model predictive control~(MPC): this is a control technique widely used to control systems under state and input constraints ~\cite{kouvaritakis2016model,RawlingsMayneDiehl2017}. To this end, the current state of the system is measured in order to predict and optimize the future system behaviour on a given (finite) prediction horizon. The optimal solution on the first portion of the considered time interval is then implemented as a control input at the plant before the whole process is repeated after obtaining a new state measurement.

While there are innumerable results regarding the stability and robustness of~MPC schemes for ordinary differential equations(ODEs), few results are known for systems governed by DAEs. The main challenge for this class of systems is the fact that input and state cannot be considered separately, so approaches known from~ODE systems based on calculating a stabilizing state feedback do not work in general.

There has been a lot of research in the related field of optimal control for~DAEs, both in an analytical context using Riccati ~\cite{Cobb83,KunkMehr08,LamoMarz13,CampKunk12} or Lur'e quations~\cite{bankmann2016linear,ReisVoig18} as well as in a numerical context~\cite{gerdts2011optimal}. However, the analytical results do not encompass state or input constraints and to our knowledge, none of these results has been extended yet to the stability analysis of~MPC schemes. The latter requires either the construction of stabilizing terminal constraints and costs or additional controllability assumptions like cost controllability, see, e.g.~\citeA{CoroGrun20} and the references therein. 

Most approaches to~MPC for~DAEs do not explicitly exploit the structure of the~DAE, but treat it as an additional constraint of the optimal control problem (OCP)~\cite{DiehBock02,bock2007constrained}. Exploiting the structure of the~DAE before discretization is done by~\citeA{YoncFind04,SjobFind07,ilchmann2019optimal}. However, these results are only applicable to regular~DAEs, and in case of~\citeA{YoncFind04,SjobFind07} some further controllability assumptions are imposed.

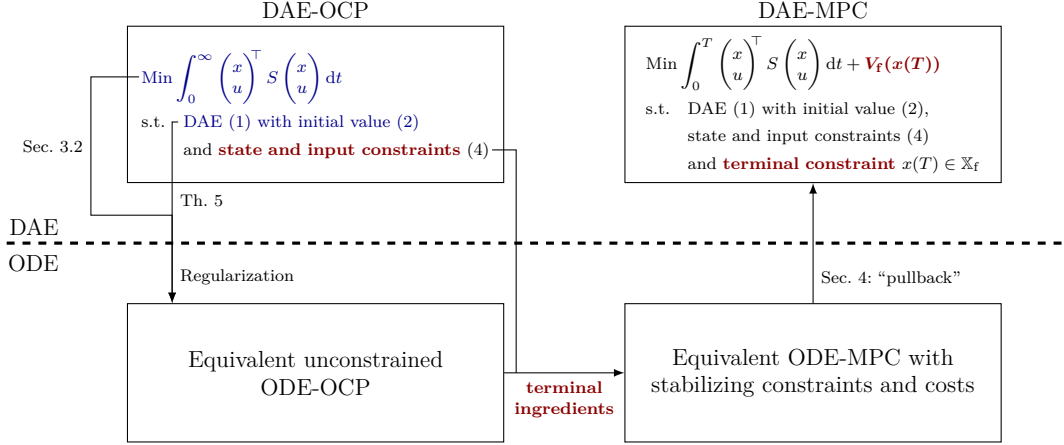
\begin{figure*}
\makeatletter
  \Hy@colorlinksfalse
  \Hy@ocgcolorlinksfalse
  \Hy@frenchlinksfalse
  \def\Hy@colorlink#1{\begingroup}%
  \def\Hy@endcolorlink{\endgroup}%
  \def\@pdfborder{0 0 0}%
  \let\@pdfborderstyle\ltx@empty
\makeatother
\centering
\begin{tikzpicture}[scale=0.8,transform shape]
\node[draw,text width=6cm,align=center,minimum height=2.6cm,font=\footnotesize] (dae-ocp)
{
	\begin{align*}
		\SwapAboveDisplaySkip
		&\color{DarkBlue}\mathrm{Min} \int_0^\infty \begin{pmatrix} x \\ u \end{pmatrix}^{\!\!\!\top} S
		\begin{pmatrix} x \\ u \end{pmatrix} \mathrm{d} t \\
		&\text{s.t.\quad \color{DarkBlue} DAE~\eqref{eq:dae} with initial value \eqref{eq:initialvalue}} \\
		&\text{\hphantom{s.t.}\quad and {\bfseries\color{DarkRed} state and input constraints}~\eqref{eq:constraints}}
	\end{align*}
};
\node[above=0cm of dae-ocp] {DAE-OCP};
\path (dae-ocp.north west) +(2mm,-8.5mm) coordinate (min) +(8.5mm,-16mm) coordinate (dae);
\path (dae-ocp.south east) +(-2mm,5.5mm) coordinate (constraints);
\path (dae-ocp.north west) +(12.5mm,-2mm) coordinate (infty);
\node[draw,below=2cm of dae-ocp,text width=6cm,align=center,minimum height=2.3cm] (ode-ocp) {Equivalent \mbox{unconstrained} ODE-OCP};
\node[draw,right=2cm of ode-ocp,text width=6cm,align=center,minimum height=2.3cm] (terminalingredients) {Equivalent ODE-MPC with stabilizing constraints and costs};
\node[draw,above=2cm of terminalingredients,text width=6cm,align=center,minimum height=2.6cm,font=\footnotesize] (mpc)
{
		\begin{align*}
		\SwapAboveDisplaySkip
		&\mathrm{Min} \int_0^T \begin{pmatrix} x \\ u \end{pmatrix}^{\!\!\!\top} S
		\begin{pmatrix} x \\ u \end{pmatrix} \mathrm{d} t + \color{DarkRed}\bm{V_\mathrm{f}(x(T))} \\
		&\text{s.t.\quad DAE \eqref{eq:dae} with initial value \eqref{eq:initialvalue}}, \\
		&\text{\hphantom{s.t.}\quad state and input constraints~\eqref{eq:constraints}} \\
		&\text{\hphantom{s.t.}\quad and {\bfseries\color{DarkRed} terminal constraint} } x(T)\in\mathbb{X}_\mathrm{f}
	\end{align*}
};
\node[above=0cm of mpc] {DAE-MPC};

\path (dae) |- coordinate (end) (ode-ocp.north);
\draw[-latex] (dae) -| node[right,pos=0.715,font=\footnotesize] {Th.~\ref{thm:odetransformation}} ($(end)-(1mm,0)$);
\draw[-latex] (min) -| node[left,pos=0.75,font=\footnotesize] {Sec.~\ref{sec:UnconstrainedOCP}} ++(-8mm,-23mm) -| node[pos=0.85,right,font=\footnotesize] {Regularization}  ($(end)-(1mm,0)$);
\draw[-latex] (ode-ocp) -- node[below,text width=2cm,align=center,font=\footnotesize\bfseries,color=DarkRed] {terminal ingredients} (terminalingredients);
\draw (constraints) -| ($(ode-ocp.east)+(2mm,0)$);
\draw[-latex] (terminalingredients) -- node[pos=0.2,right,font=\footnotesize] {Sec.~\ref{sec:stability}: \enquote{pullback}} (mpc);

\draw[dashed,very thick] (ode-ocp.north west) ++(-20mm,10mm) -- ($(terminalingredients.north east)+(10mm,10mm)$);
\node[above left=10mm and 10mm of ode-ocp.north west] {DAE};
\node[above left=4mm and 10mm of ode-ocp.north west] {ODE};
\end{tikzpicture}
\caption{Scheme of the proposed approach for~MPC for~DAEs followed in this paper}\label{fig:structure}
\end{figure*}

In the present paper, we consider~MPC for arbitrary linear time-invariant~DAEs without additional regularity or controllability assumptions. We follow the scheme presented in~\Cref{fig:structure}: in order to be able to exploit the structure of the~DAE-OCP and to use well-known results from the~MPC theory, we reduce the~DAE-OCP to an~ODE-OCP by using numerically advantageous regularization techniques and transforming the cost functional of the nominal~DAE-OCP.

On the basis of this equivalent~ODE-OCP, we construct stabilizing terminal constraints and costs for the~ODE-OCP. These terminal constraints and costs can then be transformed into terminal constraints and costs for the original~DAE-OCP. While the preceding steps of regularization and construction of the terminal ingredients for the~ODE-OCP are only adapted by us, the last step is a novelty. Its advantage is that it allows to use numerically feasible schemes to solve the~DAE-OCP, while retaining the stability guarantees which are offered by the appropriate construction of the terminal ingredients.

The paper is structured as follows: In~\Cref{sec:problemformulation}, we define the problem and propose a~MPC scheme as a solution. In~\Cref{sec:terminalingredients}, we construct suitable terminal constraints for the~MPC scheme. This is achieved by regularizing the~DAE as described in~\Cref{thm:odetransformation}, which allows to derive an equivalent~ODE-OCP in~\Cref{sec:UnconstrainedOCP}. Our main contribution is the new structural approach and the results in~\Cref{sec:stability}; this allows to prove asymptotic stability of the~MPC scheme for the nominal~DAE. We conclude with an illustrative example in~\Cref{sec:example}.

\textbf{Notation}: \(\mathcal{L}^1_\mathrm{loc}(I,\mathbb{R}^{p})\), $p \in \mathbb{N}$, denotes the space of Lebesgue-measurable functions defined on the interval \(I \subseteq \mathbb{R}\) that are locally absolutely integrable. In this context, we use the abbreviations \(\mathrm{ae}\) for \textit{almost everywhere}, and \(\mathrm{aa}\) for \textit{almost all}.
\(\mathcal{W}^{1,1}_\mathrm{loc}(I,\mathbb{R}^{p})\) denotes the Sobolev space of  weakly differentiable functions $f:I\to \mathbb{R}^{p}$ such that
\(f,\dot f \in \mathcal{L}^1_\mathrm{loc}(I,\mathbb{R}^{p})\). \(GL_n(\mathbb{R})\subset\mathbb{R}^{n\times n}\) denotes the space of invertible real matrices. For \(M\in\mathbb{R}^{n\times n}\), \(M>0\) (\(M\geq0\)) means that \(M\) is positive (semi-)definite. \(\mathds{1}\) is the all-ones vector. For \(x\in\mathbb{R}^n\), \(x\leq\mathds{1}\) means that \(x_i\leq1\) holds for all \(i\in\{1,\dots,n\}\).

\section{Problem formulation: MPC for DAEs}\label{sec:problemformulation}

\noindent We consider the differential-algebraic equation system~\([E,A,B]\in\Sigma:=\mathbb{R}^{\ell\times n}\times\mathbb{R}^{\ell\times n}\times\mathbb{R}^{\ell\times m}\) associated to
\begin{align}
	\tfrac{\mathrm{d}}{\mathrm{d}t} (E x(t)) = A x(t) + B u(t) \label{eq:dae}.
\end{align}
The system~\([E,A,B]\) is called \emph{regular} if, and only if, the linear matrix pencil~\(s E-A\) is regular, i.e. \(\ell=n\) and there exists \(\lambda\in\mathbb{C}\) such that \(\det(\lambda E-A) \neq 0\); otherwise, the system~\([E,A,B]\) is called \emph{singular}. %
The \emph{behaviour}
\begin{equation*}
	\mathfrak{B}_{[E,A,B]} := \Bigg\{ (x,u) \in \mathcal{L}_{\mathrm{loc}}^{1}(\mathbb{R}_{\geq 0},\mathbb{R}^{n+m}) \,\Bigg\vert\, \begin{array}{l} E x \in \mathcal{W}^{1,1}_{\mathrm{loc}}(\mathbb{R}_{\geq0},\mathbb{R}^\ell), \\ \text{\eqref{eq:dae} holds for almost all $t \geq 0$} \end{array} \Bigg\}
\end{equation*}
of the system~\([E,A,B]\in\Sigma\) is the set of all solutions of~\eqref{eq:dae}.
An initial value \(x^0\in\mathbb{R}^n\) is called \emph{weakly consistent} if, and only if, a solution \((x,u)\in\mathfrak{B}_{[E,A,B]}\) exists with
\begin{equation}
	(E x)(0)=E x^0. \label{eq:initialvalue}
\end{equation}
The space of \emph{weakly consistent initial values} is denoted by
\begin{equation*}
	\cIV\subseteq\mathbb{R}^n.
\end{equation*}
For every pencil~\(s E-A\in\mathbb{R}[s]^{\ell\times n}\), there exist \(T\in GL_\ell(\mathbb{R})\), \(U\in GL_n(\mathbb{R})\) that transform it into a standard form, the so-called \emph{quasi Kronecker form}~\cite{BergTren12}:
\begin{equation}
	T (s E-A) U= \begin{bmatrix} s E_\mathrm{U}-A_\mathrm{U} & 0 & 0 & 0 \\ 0 & s E_\mathrm{J}-A_\mathrm{J} & 0 & 0 \\ 0 & 0 & s E_\mathrm{N}-A_\mathrm{N} & 0 \\ 0 & 0 & 0 & s E_\mathrm{O}-A_\mathrm{O} \end{bmatrix}, \label{eq:qkf}
\end{equation}
where
\begin{itemize}
\item \(s E_\mathrm{U}-A_\mathrm{U}\in\mathbb{R}[s]^{\ell_\mathrm{U}\times n_\mathrm{U}}\), \(0\leq\ell_\mathrm{U}<n_\mathrm{U}\), \(\forall\,\lambda\in\mathbb{C}: \rank (\lambda E_\mathrm{U}-A_\mathrm{U})=\ell_\mathrm{U}\),
\item \(s E_\mathrm{J}-A_\mathrm{J}\in\mathbb{R}[s]^{n_\mathrm{J}\times n_\mathrm{J}}\), \(\rank(E_\mathrm{J})=n_\mathrm{J}\),
\item \(s E_\mathrm{N}-A_\mathrm{N}\in\mathbb{R}[s]^{n_\mathrm{N}\times n_\mathrm{N}}\), \(\forall\,\lambda\in\mathbb{C}: \rank (\lambda E_\mathrm{N}-A_\mathrm{N})=n_\mathrm{N}\), \(E_\mathrm{N}\) nilpotent,
\item \(s E_\mathrm{O}-A_\mathrm{O}\in\mathbb{R}[s]^{\ell_\mathrm{O}\times n_\mathrm{O}}\), \(\ell_\mathrm{O}>n_\mathrm{O}\geq0\), \(\forall\,\lambda\in\mathbb{C}: \rank (\lambda E_\mathrm{O}-A_\mathrm{O})=n_\mathrm{O}\).
\end{itemize}
The block sizes \(\ell_\mathrm{U}\), \(n_\mathrm{U}\), \(n_\mathrm{J}\), \(n_\mathrm{N}\), \(\ell_\mathrm{O}\), and \(n_\mathrm{O}\) are uniquely determined.

The \emph{index} of a pencil~\(s E-A\) is given by
\begin{equation*}
	\operatorname{ind} (s E-A) := \operatorname{nil\,ind} (E_\mathrm{N}) := \min\{i\in\mathbb{N} \mid E_\mathrm{N}^i=0 \}.
\end{equation*}
It can be shown that the index does not depend on the choice of transformation into quasi Kronecker form. We also define the index of the system~\([E,A,B]\in\Sigma\) as the index of the corresponding pencil~\(s E-A\).

Let the \emph{mixed state and control constraints} be given by
\begin{equation}
	\begin{bmatrix} F & G \end{bmatrix} \begin{pmatrix} x(t) \\ u(t) \end{pmatrix} \overset{\mathrm{ae}}{\leq} \mathds{1} \label{eq:constraints}
\end{equation}
with \(F\in\mathbb{R}^{p\times n}\) and~\(G\in\mathbb{R}^{p\times m}\). The constraints are only required to be fulfilled almost everywhere because~\(x\) need not be continuous. If~\(F=F' E\) holds for some~\(F'\in\mathbb{R}^{p\times \ell}\), then the constraints hold everywhere since~\(E x\) is continuous by definition of~\(\mathfrak{B}_{[E,A,B]}\).

\subsection{Singular DAEs do not allow stabilization by state feedback}\label{sec:singlexample}

\noindent Consider the singular system~\eqref{eq:dae} given by \(E=(0,1)\), \(A=(1,0)\), and \(B=1\). Its behaviour is given by
\begin{equation*}
	\mathfrak{B}_{[E,A,B]} = \left\{ \begin{pmatrix} x_1 \\ x_2 \\ u \end{pmatrix}\in\mathcal{L}^1_\mathrm{loc}(\mathbb{R}_{\geq0},\mathbb{R}^3) \,\middle|\, \dot{x}_2\overset{\text{ae}}{=}x_1+u \right\}.
\end{equation*}
Here, \(x_1\) is a free variable that can be potentially unbounded. It cannot be influenced by the control input~\(u\), which shows that for singular~DAEs, prescribing~\(u\) is generally not sufficient to achieve convergence \(\lim_{t\to\infty} x(t)=0\) of the state: for all choices \(u\in\mathcal{L}^1_\mathrm{loc}(\mathbb{R}_{\geq0},\mathbb{R})\), there exist solutions \((x,u)\in\mathfrak{B}_{[E,A,B]}\) with unbounded \(x\), e.g.
\begin{equation*}
	x_1(t)=\mathrm{e}^t, \quad x_2(t)=\int_0^t x_1(\tau)+u(\tau)\,\mathrm{d} \tau, \quad t\geq0.
\end{equation*}
Hence, in general it is necessary to have control over both the input~\(u\) as well as the (free variable part of the) state~\(x\) for singular~DAEs to achieve stabilization of the state.
Moreover, in contrast to~ODEs, it is generally impossible to construct a stabilizing linear state feedback: Let \(u = k x\) for some arbitrary but fixed \(k \in \mathbb{R}^{1 \times 2}\). Then the closed-loop system has the form
\begin{equation*}
	 \frac{\mathrm{d}}{\mathrm{d} t} \begin{pmatrix} 0, & 1 \end{pmatrix} x(t) = \begin{pmatrix} 1+k_1, & k_2 \end{pmatrix} x(t),
\end{equation*}
which is still singular. We obtain that
\begin{equation*}
	\mathfrak{B}_{[(0,1),(1+k_1,k_2),0]} = \biggl\{ \begin{psmallmatrix} x_1 \\ x_2 \end{psmallmatrix}\in\mathcal{L}^1_\mathrm{loc}(\mathbb{R}_{\geq0},\mathbb{R}^2) \,\biggm|\, \begin{array}{l} x_1 \text{ arbitrary}, \\ \dot{x}_2(t)\overset{\mathrm{ae}}{=}k_2 x_2(t)+(1+k_1) x_1(t) \end{array} \biggr\}.
\end{equation*}
Therefore, the closed-loop system admits unbounded solutions no matter how the state feedback is chosen.

\subsection{Model predictive control}\label{sec:mpc}

\noindent Our goal is to construct a feedback law such that the origin is asymptotically stable w.r.t.\,the closed-loop system while validity of the constraints~\eqref{eq:constraints} is maintained. To this end, we employ the following~MPC scheme. In every step, we measure the current state~\(E x\) of the system~\eqref{eq:dae} and solve a quadratic optimal control problem on the optimization/prediction horizon~\(T>0\). The~OCP is constrained by the system~\eqref{eq:dae} with the current system state as an initial value~\eqref{eq:initialvalue} and the mixed state and control constraints~\eqref{eq:constraints}. Clearly, solving the~OCP on an infinite-time horizon, i.e.\ $T = \infty$, would be desirable. However, due to the (mixed) state and control constraints, this is, in general, computationally intractable.

We obtain a predicted optimal solution \((\bar{x}_0^*(\cdot), \bar{u}_0^*(\cdot))|_{[0,T]}\), indicated by $\bar{\cdot}$, of which we implement only the first piece \((\bar{x}_0^*(\cdot), \bar{u}_0^*(\cdot))|_{[0,\delta)}\) up to a \emph{time shift}~\(\delta\in(0,T)\) at the plant. Note that the resulting closed-loop solution only coincides with the predicted one during the first $\delta$ time units but may differ afterwards.

After the time~\(\delta\) has passed, we repeat the procedure with the new system state.

This scheme alone would not necessarily allow for a system that is asymptotically stable w.r.t.\ the origin, see~\citeA[Sec.~1.3.4]{RawlingsMayneDiehl2017}. To guarantee stability, we incorporate additional \emph{terminal constraints}~\(\mathbb{X}_\mathrm{f} \subseteq \cIV\) and \emph{terminal costs}~\(V_\mathrm{f}:\mathbb{X}_f \rightarrow \mathbb{R}_{\geq 0}\) into the basic~MPC scheme to obtain the following algorithm. These additional degrees of freedom guarantee asymptotic stability if chosen suitably, as we will show in this paper.

\begin{algorithm*}\leavevmode\\
\textbf{Input parameters:} \([E,A,B]\in\Sigma\), \ \(F\in\mathbb{R}^{p\times n}\), \(G\in\mathbb{R}^{p\times m}\), \ \(T>0\), \ \(S\in\mathbb{R}^{(n+m)\times(n+m)}\), \\
\hphantom{\textbf{Input parameters:}} \(\mathbb{X}_\mathrm{f}\subseteq\cIV\), \(V_\mathrm{f}: \mathbb{X}_\mathrm{f}\to\mathbb{R}_{\geq0}\). \\
Set \(k=0\).
\begin{enumerate}[label={Step \arabic*:},ref={Step \arabic*},leftmargin=*]
	\item\label{item:step1} Measure \(\hat{x}_k := (Ex)(k \delta)\).
	\item\label{item:step2} Minimize \(\int_0^T \begin{psmallmatrix} \bar{x}_k(s) \\ \bar{u}_k(s) \end{psmallmatrix}^{\!\!\!\top} S \, %
		\begin{psmallmatrix} \bar{x}_k(s) \\ \bar{u}_k(s) \end{psmallmatrix}\,\mathrm{d} s + V_\mathrm{f}(\bar{x}_k(T))\) s.t.\\
	\begin{itemize}
		\item \(\frac {\mathrm{d}}{\mathrm{d}t}(E \bar{x}_k)(t) = A \bar{x}_k(t) + B \bar{u}_k(t)\)\\[-1ex]
		\item \(F \bar{x}_k(s) + G \bar{u}_k(s) \leq \mathds{1}\) for almost all \(s \in [0,T]\)\\[-1ex]
		\item \(\bar{x}_k(T) \in \mathbb{X}_\mathrm{f}\), \((E\bar{x}_k)(0) = \hat{x}_k\).
	\end{itemize}
	\item Implement \emph{first} piece \((\bar{x}_k^*(\cdot), \bar{u}_k^*(\cdot))|_{[0,\delta)}\) of optimal solution for system \([E,A,B]\) to obtain \((x(\cdot),u(\cdot))|_{[k\delta,(k+1)\delta]}\), set \(k:=k+1\), go to \ref{item:step1}.\\
\end{enumerate}
\end{algorithm*}
To construct suitable terminal constraints and terminal costs, we will consider the optimal control problem from \ref{item:step2} of the algorithm with infinite optimization horizon \(T=\infty\) and without the mixed state and control constraints~\eqref{eq:constraints}. Using regularization techniques, we will transform the~DAE constraint~\eqref{eq:dae} into an equivalent~ODE constraint so that we can solve the~OCP by well-known Riccati theory. The optimal value of this~OCP will then serve as the terminal cost \(V_\mathrm{f}: \mathbb{X}_\mathrm{f}\to\mathbb{R}_{\geq0}\), which fulfils the decrease condition
\begin{equation}\label{eq:DecentCondition}
	\forall \delta>0\,\forall \hat{x}\in\cIV: V_\mathrm{f}(x(\delta)) \leq V_\mathrm{f}(\hat{x}) - \int_0^\delta \begin{pmatrix} x(t) \\ u(t) \end{pmatrix}^{\!\!\!\top} S \begin{pmatrix} x(t) \\ u(t) \end{pmatrix} \mathrm{d} t
\end{equation}
by virtue of the Belman equation. To construct a controlled forward invariant terminal region~\(\mathbb{X}_\mathrm{f} \subseteq \cIV\), a sub-level set of~\(V_\mathrm{f}\) where the state and control constraints constraints~\eqref{eq:constraints} are fulfilled can be chosen. Together with the decrease condition~\eqref{eq:DecentCondition}, asymptotic stability of the origin w.r.t.\ the~MPC closed loop can be shown analogously to the \textit{classical}~ODE~case, see, e.g.\,\citeA{Mayne2000}.

\section{Construction of terminal ingredients}\label{sec:terminalingredients}
\noindent To construct the terminal region and terminal costs, we reduce the~DAE to an equivalent~ODE by reducing its index to~1. In order to do so, we transform the~DAE using the methods explained in the next subsection so that it fulfils the following algebraic characterization.

\begin{proposition}[{\citeA[Eq.~(3.4)]{BergReis13a}}] \label{def:index1}
The system~\([E,A,B]\in\Sigma\) has \emph{index at most one} if, and only if,
\begin{equation*}
	\image A \subseteq \image E+A \ker E.
\end{equation*}
\end{proposition}

If the~DAE has index at most one, it can be transformed to an~ODE with a state transformation.
\begin{proposition}[{\citeA[Thm.~8.1]{BennLoss15}}]\label{thm:index1}
The system~\([E,A,B]\in\Sigma\) is regular with index at most one if, and only if, there are transformation matrices \(S_\mathrm{r},T_\mathrm{r}\in GL_n(\mathbb{R})\) such that
	\begin{equation}\label{thm:svdtransformation}
		S_\mathrm{r} E T_\mathrm{r} = \begin{bmatrix} I_{\hat{n}} & 0 \\ 0 & 0 \end{bmatrix}, \quad %
		S_\mathrm{r} A T_\mathrm{r} = \begin{bmatrix} A_{11} & A_{12} \\ A_{21} & A_{22} \end{bmatrix}, \quad %
		S_\mathrm{r} B = \begin{bmatrix} B_1 \\ B_2 \end{bmatrix},
	\end{equation}
with \(\hat{n}\leq n\), \(A_{22}\in GL_{n-\hat{n}}(\mathbb{R})\).
\end{proposition}

A remark to~\Cref{thm:index1} is warranted.
\begin{remark}
\(E\) can be transformed to the form as in~\eqref{thm:svdtransformation} using a singular value decomposition: Choose orthogonal matrices $U \in \mathbb{R}^{n \times n}$ and $V \in \mathbb{R}^{n \times n}$ such that $E = U \diag(\sigma_1,\dots,\sigma_{\hat{n}},0,\dots,0) V^\top$ for \(\sigma_1,\dots,\sigma_{\hat{n}}>0\). Then \(S_\mathrm{r} := U^\top\), \(T_\mathrm{r}:= V \diag(\sigma_1^{-1},\dots,\sigma_{\hat{n}}^{-1},1,\dots,1)\) leads to~\eqref{thm:svdtransformation}.
\end{remark}

If~\([E,A,B] \in \Sigma\) is regular and has index at most one,  then the~DAE~\eqref{eq:dae} can be transformed into an~ODE with an explicit representation of the remaining states. This is made precise in the following lemma, which is an immediate consequence of~\Cref{thm:index1}.
\begin{lemma}\label{thm:odeequivalence}
	Let~\([E,A,B]\in\Sigma\) be regular with index at most one and consider the transformation into~\eqref{thm:svdtransformation}. Then it holds that
	\begin{equation}\nonumber
		\begin{psmallmatrix} x \\ u \end{psmallmatrix} \in \mathfrak{B}_{[E,A,B]}
\qquad  \Longleftrightarrow   \qquad%
		x = T_\mathrm{r} z
	\end{equation}
	with \(z = (z_1, z_2)\), where \(z_2=-A_{22}^{-1} A_{21} z_1-A_{22}^{-1} B_2 u\) and \(z_1\) solves
	\begin{equation}\nonumber
		\dot{z}_1 = (A_{11}-A_{12} A_{22}^{-1}A_{21}) z_1+(B_1-A_{12} A_{22}^{-1} B_2) u
	\end{equation}
for some $ u \in  \mathcal{L}^1_\mathrm{loc}( \mathbb{R}_{\ge 0},\mathbb{R}^m)$.
\end{lemma}

\Cref{thm:odeequivalence} essentially allows to reformulate the~DAE-OCP as an~ODE-OCP, provided that the~DAE system is regular and has index at most one. Using the regularization techniques shown in the following theorem allows to transform every~DAE system~\eqref{eq:dae} into the form form~\eqref{thm:svdtransformation}.
\begin{theorem}\label{thm:odetransformation}
Consider \([E,A,B]\in\Sigma\). Then there exist~\(\widehat{T}\in GL_{n+m}(\mathbb{R})\) and unimodular~\(U(s)=s U_1+U_0\in\mathbb{R}[s]^{\ell\times\ell}\), \(U_0,U_1\in\mathbb{R}^{\ell\times\ell}\),  such that
	\begin{equation*}
		[s E-A,-B] \, \widehat{T} = U(s) \begin{bmatrix} 0 & 0 & 0 \\ s I_{\hat{n}}-A_{11} & -A_{12} & -B_1 \\ -A_{22} & -A_{22} & -B_2 \end{bmatrix},
	\end{equation*}
	where \(\hat{n}\leq n\), \(A_{22}\in GL_{n-\hat{n}}(\mathbb{R})\).

	If~\([E,A,B]\) additionally fulfils~\eqref{eq:feedbackreg}, then \(\widehat{T}\) and \(U(s)\) can be chosen as
	\begin{equation*}
		\widehat{T} =  \begin{bmatrix} T_\mathrm{r} & 0 \\ K T_\mathrm{r} & I_m \end{bmatrix}, \quad U(s)=S_\mathrm{r}^{-1}
	\end{equation*}
	for some \(K\in\mathbb{R}^{m\times n}\), \(S_\mathrm{r},T_\mathrm{r}\in GL_n(\mathbb{R})\).

\begin{proof}
This will be proved in the following~\Cref{sec:regularizationproof}: The first part of the assertion follows immediately from~\Cref{thm:unimodularregularization} in combination with~\Cref{thm:index1}, while the second part is a direct consequence of~\Cref{thm:feedbackregularization} and~\Cref{thm:index1}.
\end{proof}
\end{theorem}

\subsection{Proof of \Cref{thm:odetransformation}}\label{sec:regularizationproof}

\noindent We use two different approaches to regularize the system~\eqref{eq:dae}: If the system fulfils a certain rank condition, we use the method by~\citeA{BunsMehr92,BunsMehr94} to find a regularizing feedback. If it is not, we use the more general approach proposed by~\citeA{BergDoor15}. While the latter approach is also applicable to the case where the rank condition is fulfilled, the former is in this case numerically more attractive.

\begin{definition}[{\citeA[Def.~2.1]{BergReis13a}}]
\([E,A,B]\in\Sigma\) is called \emph{impulse controllable} if, and only if,
\begin{equation*}
	\forall\, x^0\in\mathbb{R}^n\,\exists\, (x,u)\in\mathfrak{B}_{[E,A,B]}: E x^0=E x(0).
\end{equation*}
\end{definition}

\begin{proposition}[{\citeA[Rem.~4.2]{BergReis13a}}]
If \([E,A,B]\in\Sigma\) is regular, then it is impulse controllable if, and only if,
\begin{equation}
	\rank [E,A Z,B]=n \quad\text{for some}\quad Z\in\mathbb{R}^{n\times(n-\rank E)}, \  \image Z=\ker E.\label{eq:feedbackreg}
\end{equation}
\end{proposition}

Next, we present a well-known technique to regularize the system~\eqref{eq:dae}.
\begin{proposition}\label{thm:feedbackregularization}
	\([E,A,B]\in\Sigma\) fulfils~\eqref{eq:feedbackreg} if, and only if, there exists a feedback matrix~\(K \in \mathbb{R}^{n \times m}\) %
	such that~\([E,A+B K,B]\) is regular and has index at most one. Moreover, the behaviours are linked by
	\begin{equation*}
		\begin{pmatrix} x \\ u \end{pmatrix}\in\mathfrak{B}_{[E,A,B]} \ \Longleftrightarrow\ \begin{pmatrix} x \\ v \end{pmatrix}\in\mathfrak{B}_{[E,A+B K,B]},\ \text{where } u=K x+v.
	\end{equation*}
\begin{proof}
	Sufficiency of the first assertion is proved in~\citeA[Theorem~6]{BunsMehr92}, while necessity is shown in~\citeA[Theorem~4]{BunsMehr94}. %
	The second statement is immediate.
\end{proof}
\end{proposition}
The behaviour of the nominal system~\eqref{eq:dae} and the regularized system are coupled by an input transformation.

If \([E,A,B]\in\Sigma\) does not fulfill~\eqref{eq:feedbackreg}, 
we will conduct the regularization by special unimodular transformations.
A polynomial matrix~\(U(s)\in\mathbb{R}[s]^{n\times n}\) is called \emph{unimodular} if, and only if, \(U^{-1}(s)\in\mathbb{R}[s]^{n\times n}\) exists
such that $U(s) \, U^{-1}(s)=I_n$.
Clearly, \(U(s)\) is unimodular if, and only if, its determinant is a constant nonzero polynomial.

In the next proposition, we show that the nominal system~\eqref{eq:dae} can be regularized by an input transformation.
\begin{proposition}\label{thm:unimodularregularization}
	For any system~\([E,A,B]\in\Sigma\) there exist~\(\widehat{T}\in GL_{n+m}(\mathbb{R})\) and %
	unimodular~\(U(s)=s U_1+U_0\in\mathbb{R}[s]^{\ell\times\ell}\), \(U_0,U_1\in\mathbb{R}^{\ell\times\ell}\),  such that
	\begin{equation}
		[s E-A,-B] \, \widehat{T} = U(s) \begin{bmatrix} 0 & 0 \\ s E_\mathrm{r}-A_\mathrm{r} & -B_\mathrm{r} \end{bmatrix}, 
\label{eq:bergdoortrans}
	\end{equation}
	where \(s E_\mathrm{r}-A_\mathrm{r}\in\mathbb{R}[s]^{r\times r}\), \(r\leq\ell\) is regular and has index at most one, \(B_\mathrm{r}\in\mathbb{R}^{\ell\times(n+m-r)}\). Then, the following implications hold
	\begin{enumerate}[label={(\roman*)}]
		\item\label{item:behaviourregularized} \(\begin{psmallmatrix} x \\  u \end{psmallmatrix}\in\mathfrak{B}_{[E,A,B]}\) \qquad\(\Longrightarrow\)\qquad \(\widehat{T}^{-1} \begin{psmallmatrix} x \\  u \end{psmallmatrix}\in\mathfrak{B}_{[E_\mathrm{r},A_\mathrm{r},B_\mathrm{r}]}\). \\[-3mm]
		\item\label{item:smoothnessregularized} If \(\begin{psmallmatrix} x \\ u \end{psmallmatrix} \in \mathfrak{B}_{[E_\mathrm{r},A_\mathrm{r},B_\mathrm{r}]}\) satisfies 
			\begin{equation}\nonumber
				 \Bigl(U_0 \begin{bsmallmatrix} 0 & 0 \\ E_\mathrm{r} & 0 \end{bsmallmatrix} - U_1 \begin{bsmallmatrix} 0 & 0 \\ A_\mathrm{r} & B_\mathrm{r} \end{bsmallmatrix}\Bigr) \begin{psmallmatrix} x \\ u \end{psmallmatrix} \in\mathcal{W}^{1,1}_\mathrm{loc},
			\end{equation}
			then \(\widehat{T} \begin{psmallmatrix} x \\ u \end{psmallmatrix}\in\mathfrak{B}_{[E,A,B]}\) holds.
	\end{enumerate}
\begin{proof}
	The first part of the proposition  is shown in~\citeA{BergDoor15}. Note that it also follows from~\citeA{BergDoor15} that \(U(s)\), \(E_\mathrm{r}\) and \(A_\mathrm{r}\) are chosen such that no quadratic term occurs on the right-hand side of~\eqref{eq:bergdoortrans}. Assertion~\ref{item:smoothnessregularized} can be immediately concluded from
	\begin{align*}
		\mathrel{\hphantom{=}} [s E-A,-B] \, \widehat{T} = s %
		\Bigl(U_0 \begin{bsmallmatrix} 0 & 0 \\ E_\mathrm{r} & 0 \end{bsmallmatrix}-U_1 \begin{bsmallmatrix} 0 & 0 \\ A_\mathrm{r} & B_\mathrm{r} \end{bsmallmatrix}\Bigr) %
		-U_0 \begin{bsmallmatrix} 0 & 0 \\ A_\mathrm{r} & B_\mathrm{r} \end{bsmallmatrix}.
	\end{align*}
	
	To show Assertion~\ref{item:behaviourregularized}, let \(\begin{psmallmatrix} x \\  u \end{psmallmatrix}\in\mathfrak{B}_{[E,A,B]}\) be arbitrary, %
	so by definition \(E x\in\mathcal{W}^{1,1}_\mathrm{loc}(\mathbb{R}_{\geq0},\mathbb{R}^\ell)\). %
	We need to show that~\([E_\mathrm{r},0_{r\times(n+m-r)}] \widehat{T}^{-1} \begin{psmallmatrix} x \\ u \end{psmallmatrix}\in\mathcal{W}^{1,1}_\mathrm{loc}(\mathbb{R}_{\geq0},\mathbb{R}^r)\). %
	By following the proof outlined in~\citeA{BergDoor15}, we obtain
	\begin{equation*}
		\ker [E,0_{\ell\times m}] \widehat{T} \subseteq \ker \begin{bmatrix} 0 & 0 \\ E_\mathrm{r} & 0 \end{bmatrix}.
	\end{equation*}
	Hence
	\begin{align*}
		\mathrel{\hphantom{=}} \begin{bmatrix} 0 & 0 \\ E_\mathrm{r} & 0 \end{bmatrix} \widehat{T}^{-1} \begin{psmallmatrix} x \\ u \end{psmallmatrix}
		= %
		\underbrace{[E,0] \begin{psmallmatrix} x \\  u \end{psmallmatrix}}_{\in\mathcal{W}^{1,1}_\mathrm{loc}(\mathbb{R}_{\geq0},\mathbb{R}^\ell)} + %
		\underbrace{\Bigl(\begin{bmatrix} 0 & 0 \\ E_\mathrm{r} & 0 \end{bmatrix}-[E,0] \widehat{T}\Bigr) \widehat{T}^{-1} \begin{psmallmatrix} x \\  u \end{psmallmatrix}}_{=0},
	\end{align*}
	which proves Assertion~\ref{item:behaviourregularized}.
\end{proof}
\end{proposition}

The next example
illustrates~\Cref{thm:unimodularregularization}
for a regular index two~DAE.
\begin{example}
	Consider the regular~DAE
	\begin{equation}
		\frac{\mathrm{d}}{\mathrm{d} t} \begin{bmatrix} 0 & 1 \\ 0 & 0 \end{bmatrix} x(t) = \begin{bmatrix} 1 & 0 \\ 0 & 1 \end{bmatrix} x(t). \label{eq:examplenilpotent}
	\end{equation}
	Its behaviour is given by \((x,u)\in \mathcal{L}^1_\mathrm{loc}(\mathbb{R}_{\geq0},\mathbb{R}^2)\) %
	satisfying \(x_2\in\mathcal{W}^{1,1}_\mathrm{loc}(\mathbb{R}_{\geq0},\mathbb{R})\), \(x_1\overset{\mathrm{ae}}{=}0\), and \(x_2=0\). Setting
	\begin{equation*}
		\widehat{T} = I_2, \quad U(s) = \begin{bmatrix*}[r] -1 & s \\ 0 & -1 \end{bmatrix*},
	\end{equation*}
	the system~\eqref{eq:examplenilpotent} can be transformed, as described in~\Cref{thm:unimodularregularization}, into the regular index~1 system \([0_{2\times2},I_2]\). Its behaviour is
	\begin{equation*}
		\bigl\{\begin{psmallmatrix} x \\ u \end{psmallmatrix}\in\mathcal{L}^1_\mathrm{loc}(\mathbb{R}_{\geq0},\mathbb{R}^2) \,\bigm|\, %
		x_1\overset{\mathrm{ae}}{=}0,\ x_2\overset{\mathrm{ae}}{=}0 \bigr\}.
	\end{equation*}
	Note that the component~\(x_2\) is required to be smoother in the nominal~DAE than in the regularized~DAE,
	as it needs to be~\(0\) everywhere instead of \emph{almost} everywhere.
\end{example}

In the following example a singular~DAE (with over- and underdetermined blocks) is considered.
\begin{example}
	Consider the singular system
	\begin{equation}\label{eq:exsingular}
		\frac{\mathrm{d}}{\mathrm{d} t} \begin{bmatrix} 0 & 1 & 0 \\ 0 & 0 & 0 \\ 0 & 0 & 1 \end{bmatrix} x(t) = %
		\begin{bmatrix} 1 & 0 & 0 \\ 0 & 0 & 1 \\ 0 & 0 & 0 \end{bmatrix} x(t) + \begin{pmatrix} 0 \\ 0 \\ 1 \end{pmatrix} u(t).
	\end{equation}
	Using~\Cref{thm:unimodularregularization}, the transformation
	\begin{equation*}
		\begin{pmatrix} y_1 \\ y_2 \\ y_3 \\ v \end{pmatrix} := \begin{bmatrix} 0 & 0 & 1 & 0 \\ 0 & 1 & 0 & 0 \\ 0 & 0 & 0 & 1 \\ 1 & 0 & 0 & 0 \end{bmatrix} \begin{pmatrix} x_1 \\ x_2 \\ x_3 \\ u \end{pmatrix}
     \end{equation*}
     yields the regular index 1 system
	\begin{equation*}
		\frac{\mathrm{d}}{\mathrm{d} t} \begin{bmatrix} 0 & 0 & 0 \\ 0 & 1 & 0 \\ 0 & 0 & 0 \end{bmatrix} y(t) = %
		\begin{bmatrix} 1 & 0 & 0 \\ 0 & 0 & 0 \\ 0 & 0 & -1 \end{bmatrix} y(t) + \begin{pmatrix} 0 \\ 1 \\ 0 \end{pmatrix} v(t).
	\end{equation*}
	Hence the equivalent~ODE is given by
	\begin{equation*}
		\dot{z}_1(t) = v(t), \quad z_{2,1}(t)=z_{2,2}(t)=0,
	\end{equation*}
	where~\(z_1=y_2\), \(z_{2,1}=y_1\), \(z_{2,2}=y_3\).
	We see that in the singular system~\eqref{eq:exsingular}, \(x_1\) can be chosen freely and is therefore more of an input than a state, which is reflected by the regularized system where \(v=x_1=z_1\). In turn, the \enquote{input}~\(u\) is~\(0\) almost everywhere, hence it is a state \(u=y_3=z_{2,2}\) of the regularized system.
\end{example}

\subsection{Optimal control for DAEs without constraints}\label{sec:UnconstrainedOCP}

\noindent We consider the optimal control problem
\begin{align}
	\mathrm{Min} \int_0^T \begin{pmatrix} x(t) \\ u(t) \end{pmatrix}^{\!\!\!\top} S
	\begin{pmatrix} x(t) \\ u(t) \end{pmatrix} \mathrm{d} t 
	\quad\text{subject to  \eqref{eq:dae}}\label{eq:ocp}
\end{align}
with \textit{optimization horizon} \(T\in\mathbb{R}\cup\{\infty\}\), 
\(S=S^\top\in\mathbb{R}^{(n+m)\times (n+m)}\), and \textit{cost functional}
\begin{align*}
	J_T:\nobreak\mathfrak{B}_{[E,A,B]} \to \mathbb{R}\cup\{\pm\infty\}, \quad (x,u) &\mapsto \int_0^T \begin{pmatrix} x(t) \\ u(t) \end{pmatrix}^{\!\!\!\top} S
	\begin{pmatrix} x(t) \\ u(t) \end{pmatrix} \mathrm{d} t.
\end{align*}
The \textit{value function} for the OCP without state and control constraints is defined by
\begin{equation}
	\begin{aligned}
	V_T: \cIV &\to \mathbb{R} \cup \{\pm\infty\} \\
	x^0 &\mapsto 
	\inf  \setdef{J_T(x,u) }{  (x,u)\in\mathfrak{B}_{[E,A,B]} \  \text{with}  \  (E x)(0) = E x^0}.
	\end{aligned}
	\label{eq:daeoptimalvalue}
\end{equation}

In order to solve this~OCP, we transform it into an equivalent~OCP that is constrained by an~ODE instead of the nominal DAE with the help of~\Cref{thm:odetransformation}. To this end, define for \([E,A,B]\in\Sigma\),
\begin{equation}
	\widehat{A} := A_{11}-A_{12} A_{22}^{-1}A_{21}, \quad
	\widehat{B} := B_1 - A_{12} A_{22}^{-1} B_2,
	\label{eq:equivalentode}
\end{equation}
where \(A_{11}\), \(A_{12}\), \(A_{21}\), \(A_{22}\), \(B_1\) and \(B_2\) are defined as in~\Cref{thm:odetransformation}.

Then the~ODE-OCP we want to consider is given by the cost functional
\begin{equation*}
	\widehat{J}_T: \mathfrak{B}_{[I_{\hat{n}},\hat{A},\hat{B}]}\to\mathbb{R}\cup\{\pm\infty\},
	\quad
	\begin{psmallmatrix} z_1 \\ v \end{psmallmatrix} \mapsto \int_0^T \begin{pmatrix} z_1(t) \\ v(t) \end{pmatrix}^{\!\!\!\top} \widehat{S}
	\begin{pmatrix} z_1(t) \\ v(t) \end{pmatrix} \mathrm{d} t,
\end{equation*}
where, 
for~$\widehat{T}$ as in~\eqref{eq:bergdoortrans},
\begin{equation}\label{eq:X}
	\widehat{S} :=  X^{-1} S X,
	\quad
	X := \widehat{T} \begin{bmatrix}  I_{\hat{n}} & 0 \\ -A_{22}^{-1} A_{21} & -A_{22}^{-1} B_2 \\ 0  & I_m \end{bmatrix}
\end{equation}
using the notation from~\Cref{thm:odetransformation}. Its value function is given by
\begin{equation}\label{eq:odeocp}
\begin{array}{rcl}
	\widehat{V}_T: \cIV &\to & \mathbb{R}\cup\{\pm\infty\}  \\
	x^0 &\mapsto& \inf 
\setdef{J_T(z_1,v) }{  
\begin{array}{l}
(z_1,v)\in\mathfrak{B}_{[I_{\hat{n}},\hat{A},\hat{B}]} \  \text{with}  \\ 
z_1(0)\in[I_{\hat{n}},0] X^{-1}  
\begin{psmallmatrix} \{x^0\} \\ \mathbb{R}^m  \end{psmallmatrix}
\end{array}
}.
\end{array}
\end{equation}

We obtain that the value function~\eqref{eq:daeoptimalvalue} of the~DAE-OCP and its counterpart~\eqref{eq:odeocp} of the~ODE-OCP coincide, as detailed in the following theorem.
\begin{theorem}\label{thm:odeocp}
Consider \([E,A,B]\in\Sigma\) and the optimal control problems defined by~\eqref{eq:daeoptimalvalue} and~\eqref{eq:odeocp}. Then
\begin{equation*}
	V_T \equiv \widehat{V}_T.
\end{equation*}
If, in addition, one of the following conditions holds:
\begin{enumerate}[label={(\roman*)}]
\item\label{enum:cond_regular} \([E,A,B]\in\Sigma\) is regular, or
\item\label{enum:cond_nooverdet} the extended system \([[E,0],[A,B]]\) does not contain an overdetermined part, i.e. \(s E_\mathrm{O}-A_\mathrm{O}\) is void in any transformation of \([[E,0],[A,B]]\) into quasi Kronecker form~\eqref{eq:qkf},
\end{enumerate}
then,
for~$\widehat{T}$ and~$X$  as in~\eqref{eq:bergdoortrans}
and~\eqref{eq:X}, resp.,
\begin{equation}
	[I_{\hat{n}},0] X^{-1} \begin{pmatrix} \{x^0\} \\ \mathbb{R}^m \end{pmatrix} = \biggl\{[I_{\hat{n}},0] \widehat{T}^{-1} \begin{pmatrix} x^0 \\ 0 \end{pmatrix}\biggr\}, \label{eq:constraintroom}
\end{equation}
i.e. we can replace the initial constraint \(z_1(0)\in[I_{\hat{n}},0] X^{-1} \begin{psmallmatrix} \{x^0\} \\ \mathbb{R}^m \end{psmallmatrix}\) by an initial condition \(z_1(0)=[I_{\hat{n}},0] \widehat{T}^{-1} \begin{psmallmatrix} x^0 \\ 0 \end{psmallmatrix}\).

\begin{proof}
Consider \((x,u)\in\mathfrak{B}_{[E,A,B]}\). By~\Cref{thm:odetransformation} and~\Cref{thm:odeequivalence}, it follows for
\begin{equation*}
	\begin{pmatrix} z_1 \\ v \end{pmatrix} := \begin{bmatrix} I_{\hat{n}} & 0 \\ -A_{22}^{-1} A_{21} & -A_{22}^{-1} B_2 \\ 0  & I_m \end{bmatrix}^{-1} \widehat{T}^{-1} \begin{pmatrix} x \\ u \end{pmatrix} = X^{-1} \begin{pmatrix} x \\ u \end{pmatrix}
\end{equation*}
that \(\begin{psmallmatrix} z_1 \\ v \end{psmallmatrix}\in\mathfrak{B}_{[I_{\hat{n}},\hat{A},\hat{B}]}\). 
Let \(\begin{psmallmatrix} x \\ u \end{psmallmatrix}\in\mathfrak{B}_{[E,A,B]}\) be arbitrary.
Substituting this in~\eqref{eq:ocp} yields
\begin{align*}
	\begin{pmatrix} x \\ u \end{pmatrix}^{\!\!\!\top} S \begin{pmatrix} x \\ u \end{pmatrix}
	&= \begin{pmatrix} z_1 \\ v \end{pmatrix}^\top \begin{bmatrix}  I_{\hat{n}} & 0 \\ -A_{22}^{-1} A_{21} & -A_{22}^{-1} B_2 \\ 0  & I_m \end{bmatrix}^\top \widehat{T}^\top \, S \, \widehat{T} \begin{bmatrix}  I_{\hat{n}} & 0 \\ -A_{22}^{-1} A_{21} & -A_{22}^{-1} B_2 \\ 0  & I_m \end{bmatrix} \begin{pmatrix} z_1 \\ v \end{pmatrix} \\
	&= \begin{pmatrix} z_1 \\ v \end{pmatrix}^\top \widehat{S} \, \begin{pmatrix} z_1 \\ v \end{pmatrix},
\end{align*}
therefore
\begin{equation*}
	J_T(\begin{psmallmatrix} x \\ u \end{psmallmatrix}) = \widehat{J}_T(\begin{psmallmatrix} z_1 \\ v \end{psmallmatrix}).
\end{equation*}
Now following the same argument as~\citeA[Thm.~1]{ilchmann2019optimal}, we conclude that
\begin{equation*}
	\forall x^0\in\cIV: V_T(x^0)=\widehat{V}_T(x^0),
\end{equation*}
which proves the assertion.

If instead the additional condition~\ref{enum:cond_nooverdet} holds, we assume without loss of generality that \([E,A]\) is given in Kronecker form, i.e.
\begin{align*}
	s [E,0]-[A,B] &= s \left[\scalebox{0.8}{$\begin{array}{ccc|c} \diag(K_{n_1},\dots,K_{n_u}) & 0 & 0 \\ 0 & I_{n_J} & 0 & 0 \\ 0 & 0 & N \end{array}$}\right] - \left[\scalebox{0.8}{$\begin{array}{ccc|c} \diag(L_{n_1},\dots,L_{n_u}) & 0 & 0 \\ 0 & J & 0 & B \\ 0 & 0 & I_{n_N} \end{array}$}\right] \\
	&\in\mathbb{R}[s]^{\ell\times(\ell+u+m)},
\end{align*}
where
\begin{equation*}
\begin{lgathered}
s K_{n_i}-L_{n_i}=\begin{bmatrix} s & 1 & & 0 \\ & \ddots & \ddots &  \\ 0 & & s & 1 \end{bmatrix} \in \mathbb{R}^{(n_i-1)\times n_i},\ i=1,\dots,u \quad\text{and} \\
N \text{ nilpotent.}
\end{lgathered}
\end{equation*}
Writing
\begin{align*}
e_i & :=  (0_{1\times(i-1)},1,0_{1\times(n_1+\dots+n_u-i)})^\top\in\mathbb{R}^{n_1+\dots+n_u}, 
\qquad  i= 1,2,\dots,n_1+\dots+n_u\\
	E_\mathrm{r} &:= \begin{bmatrix} I_{\hat{n}} & 0 \\ 0 & 0_{n_N\times n_N} \end{bmatrix}\in\mathbb{R}^{\ell\times\ell}, 
	\qquad \hat{n}:=\ell-n_N, \\
	A_\mathrm{r} &:= \begin{bmatrix} \diag(N_{n_1-1},\dots,N_{n_u-1}) & 0 & 0 \\ 0 & J & 0 \\ 0 & 0 & I_{n_N} \end{bmatrix}\in\mathbb{R}^{\ell\times\ell}, \\
	B_\mathrm{r} &:= \left[\begin{array}{ccc|c} -e_{n_1-1} & \dots & -e_{(n_1-1)+\dots+(n_u-1)} & \multirow{2}{*}{$B$} \\ 0_{(n_J+n_N-u)\times1} & \dots & 0_{(n_J+n_N-u)\times1} \end{array}\right]\in\mathbb{R}^{\ell\times(u+m)}, \\
	\widehat{T} &:= \begin{bmatrix} \diag(K_{n_1}^\top,\dots,K_{n_u}^\top) & 0 & e_{n_1} & \dots & e_{n_1+\dots+n_u} & 0 \\ 0 & I_{n_J+n_N} & 0 & \dots & 0 & 0 \\ 0 & 0 & 0 & \dots & 0 & I_m \end{bmatrix}\in\mathbb{R}^{(\ell+u+m)\times(\ell+u+m)}, \\
	U(s) &:= \begin{bmatrix} I_{n_1+\dots+n_u-u} & 0 & 0 \\ 0 & I_{n_J} & 0 \\ 0 & 0 & -(s N-I_{n_N}) \end{bmatrix}\in\mathbb{R}[s]^{\ell\times\ell}
\end{align*}
yields
\begin{align*}
	\mathrel{\hphantom{=}} [s E-A,-B] \widehat{T}
	&= s \left[\begin{array}{ccc|c|c} I_{(n_1-1)+\dots+(n_u-1)} & 0 & 0 & & \\ 0 & I_{n_J} & 0 & 0 & 0 \\ 0 & 0 & N &  \end{array}\right] \\
	&\mathrel{\hphantom{=}} -  \left[\scalebox{0.8}{$\begin{array}{ccc|ccc|c} \diag(N_{n_1-1},\dots,N_{n_u-1}) & 0 & 0 & e_{n_1-1} & \dots & e_{(n_1-1)+\dots+(n_u-1)} \\ 0 & J & 0 & 0_{(n_J-u)\times1} & \dots & 0_{(n_J-u)\times1} & -B \\ 0 & 0 & I_{n_N} & 0_{n_N\times1} & \dots & 0_{n_N\times1} \end{array}$}\right] \\
	&= U(s) \Biggl(s \left[\begin{array}{ccc|c|c} I_{(n_1-1)+\dots+(n_u-1)} & 0 & 0 & & \\ 0 & I_{n_J} & 0 & 0 & 0 \\ 0 & 0 & 0_{n_N\times n_N} &  \end{array}\right] \\
	&\mathrel{\hphantom{=}} -  \left[\scalebox{0.8}{$\begin{array}{ccc|ccc|c} \diag(N_{n_1-1},\dots,N_{n_u-1}) & 0 & 0 & e_{n_1-1} & \dots & e_{(n_1-1)+\dots+(n_u-1)} \\ 0 & J & 0 & 0_{(n_J-u)\times1} & \dots & 0_{(n_J-u)\times1} & -B \\ 0 & 0 & I_{n_N} & 0_{n_N\times1} & \dots & 0_{n_N\times1} \end{array}$}\right]\Biggr) \\
	&= U(s) [s E_\mathrm{r}-A_\mathrm{r}, -B_\mathrm{r}].
\end{align*}
The system \([E_\mathrm{r},A_\mathrm{r},B_\mathrm{r}]\in\mathbb{R}^{\ell\times\ell}\times\mathbb{R}^{\ell\times\ell}\times\mathbb{R}^{\ell\times(m+u)}\) is regular with index~1. 
Since  \(\widehat{T}^{-1}=\widehat{T}^\top\), 
it follows that
\begin{align*}
	[I_{\hat{n}},0] X^{-1} &\overset{\eqref{eq:X}}{=} [I_{\hat{n}},0_{\hat{n}\times(u+m)}] \begin{bmatrix} I_{\hat{n}} & 0 \\ 0_{n_N\times\hat{n}} & * \\ 0 & I_{m+u} \end{bmatrix}^{-1} \widehat{T}^{-1} \\
	&= [I_{\hat{n}},0_{\hat{n}\times(u+m)}] \begin{bmatrix} I_{\hat{n}} & 0_{\hat{n}\times n_N} & 0 \\ 0 & 0_{(m+u)\times n_N} & I_{m+u} \end{bmatrix} \widehat{T}^{-1} \\
	&= \begin{bmatrix} I_{\hat{n}} & 0_{\hat{n}\times(n_N+m+u)} \end{bmatrix} \widehat{T}^\top = \begin{bmatrix} \diag(K_{n_1},\dots,K_{n_u}) & 0_{\hat{n}\times(n_J+n_N+m)} \end{bmatrix}.
\end{align*}
Therefore,  
\begin{equation*}
\forall\, x^0\in\mathbb{R}^n\ \forall\, u^0\in\mathbb{R}^m\ : \
	[I_{\hat{n}},0] X^{-1} \begin{pmatrix} x^0 \\ u^0 \end{pmatrix} = [I_{\hat{n}},0] X^{-1} \begin{pmatrix} x^0 \\ 0 \end{pmatrix},
\end{equation*}
which proves~\eqref{eq:constraintroom}.
\end{proof}
\end{theorem}

\textbf{Conjecture} \ 
We strongly believe~-- although we could not prove it~--
that~\eqref{eq:constraintroom} holds without any of the assumptions~(i)
or~(ii) in Theorem~\ref{thm:odeocp}.\\

We have transformed the~DAE-OCP into an equivalent~ODE-OCP. In order to ensure existence and uniqueness of an optimal control trajectory, i.e.
\begin{gather*}
	\forall\, x^0\in\cIV\,\exists\, \text{unique } \begin{psmallmatrix} x^* \\  u^* \end{psmallmatrix}\in\mathfrak{B}_{[E,A,B]}
	\quad  \cap  \quad \mathcal{C}^\infty(\mathbb{R}_{\geq0},\mathbb{R}^{n+m}): \\
	(E x^*)(0)=E x^0 \land J_T(x^*,u^*)=V_T(x^0),
\end{gather*}
we impose some standard assumptions on the~ODE-OCP~\cite{LancRodm95}: for the system~\([E,A,B]\in\Sigma\),  define
\begin{equation}
	\widehat{A} := A_{11} - A_{12}A_{22}^{-1}A_{21}, \quad \widehat{B} := B_1 - A_{12}A_{22}^{-1} B_2 \label{eq:abhat}
\end{equation}
where \(A_{ij}\), \(B_i\), \(i\in\{1,2\}\) are given by~\eqref{thm:svdtransformation}. Furthermore partition~\(\widehat{S}\), defined in~\eqref{eq:X}, as
\begin{equation}
	\widehat{S} := \begin{bmatrix} \widehat{Q} & \widehat{H} \\ \widehat{H}^\top & \widehat{R} \end{bmatrix}, \quad \widehat{Q}\in\mathbb{R}^{\hat{n}\times\hat{n}}. \label{eq:shat}
\end{equation}
Then we require the following assumption.
\begin{assumption}
Assume that the following properties hold for~\(\widehat{A}\), \(\widehat{B}\), \(\widehat{S}\) as defined in~\eqref{eq:abhat} and~\eqref{eq:shat}.
\begin{itemize}\label{assumption}
	\item\label{ass:shatpossemidef} \(\widehat{S}\geq0\),
	\item\label{ass:stabilizable} the pair \((\widehat{A},\widehat{B})\) is stabilizable,
	\item\label{ass:rhatpos} \(\widehat{R} = [0,I_{n-\hat{n}}] \, \widehat{S} \, [0,I_{n-\hat{n}}]^\top >0\),
	\item\label{ass:observable} \(\bigl( \widehat{A}, \widehat{Q} \bigr)\) \ is observable,
	\item\label{ass:shatnondegenerate} \(\rank \widehat{S} = \rank (\widehat{Q} + \widehat{R})\).
\end{itemize}
\end{assumption}

\begin{proposition}[{\citeA[Prop.~16.2.8]{LancRodm95}}]
Assume that~\Cref{assumption} holds, and consider for \(\widehat{A}, \widehat{B}, \widehat{Q}, \widehat{H}, \widehat{R}\) as defined in~\eqref{eq:abhat} and~\eqref{eq:shat} the algebraic Riccati equation
\begin{equation}
	\widehat{A}^\top \widehat{P}+\widehat{P} \widehat{A}+\widehat{Q}-(\widehat{P} \widehat{B}+\widehat{H}) \widehat{R}^{-1} (\widehat{P} \widehat{B}+\widehat{H})^\top = 0. \label{eq:are}
\end{equation}
Then this equation has a unique solution~\(\widehat{P}=\widehat{P}^\top\in\mathbb{R}^{\hat{n}\times\hat{n}}\) and this solution satisfies \(\widehat{P}>0\).
\end{proposition}

We obtain the desired existence and uniqueness result on the optimal control.
\begin{proposition}\label{thm:existenceoptimaltrajectory}
Consider the system~\([E,A,B]\in\Sigma\), and assume that~\Cref{assumption} holds. Then
\begin{equation*}
	\begin{lgathered}
	\forall\, x^0\in\cIV\,\exists\,\text{unique } \begin{psmallmatrix} x^* \\  u^* \end{psmallmatrix}\in\mathfrak{B}_{[E,A,B]}:	(E x^*)(0)=E x^0 \land J_T(x^*,u^*)=V_T(x^0),
	\end{lgathered}
\end{equation*}

\begin{proof}
Consider the~ODE-OCP
\begin{align}
	&\mathrm{Minimize} \int_0^T \begin{pmatrix} z_1(t) \\ v(t) \end{pmatrix}^{\!\!\!\top} \widehat{S}
	\begin{pmatrix} z_1(t) \\ v(t) \end{pmatrix} \mathrm{d} t \label{eq:odeocpProof} \\
	&\text{s.t. } \dot{z}_1(t) = \widehat{A} z_1(t) + \widehat{B} v(t), \ z_1(0) \in [I_{\hat{n}},0] \, X^{-1} \begin{psmallmatrix} \{x^0\} \\ \mathbb{R}^m \end{psmallmatrix}. \label{eq:odeivpProof}
\end{align}
For \(T<\infty\), this optimal control problem has a unique solution~\((z^*,u^*)\in\mathcal{C}^\infty(\mathbb{R}_{\geq0},\mathbb{R}^n\times\mathbb{R}^m)\) %
according to~\citeA[Theorem~16.4.2]{LancRodm95}; for \(T=\infty\), the same result follows from~\citeA[Theorem~16.3.3]{LancRodm95}. %
According to~\Cref{thm:unimodularregularization}, it follows that~\(\hat{T} \begin{psmallmatrix} z^* \\ u^* \end{psmallmatrix}\in\mathfrak{B}_{[E,A,B]}\), hence
\begin{equation*}
	V_T(x^0) = \widehat{V}_T(x^0) = \widehat{J}_T\Bigl(\begin{psmallmatrix} x^* \\  u^* \end{psmallmatrix}\Bigr) = J_T\Bigl(T \begin{psmallmatrix} z^* \\  u^* \end{psmallmatrix}\Bigr).
\end{equation*}
This proves the assertion.
\end{proof}
\end{proposition}

Furthermore, we can prove the Bellman equation for the~DAE-OCP:
\begin{proposition}\label{thm:bellman}
Consider the system~\([E,A,B]\in\Sigma\). Let \(x^0\in\cIV\) be arbitrary and \(\begin{psmallmatrix} x^* \\ u^* \end{psmallmatrix}\in\mathfrak{B}_{[E,A,B]}\) with \((E x^*)(0)=x^0\) be an optimal trajectory, i.\,e. \(J_\infty(\begin{psmallmatrix} x^* \\ u^* \end{psmallmatrix})=V_\infty(x^0)\). Then for all \(T>0\), it holds that
\begin{align}
	V_\infty(x^0) &= J_T\Bigl(\begin{psmallmatrix} x^* \\ u^* \end{psmallmatrix}\Bigr) + V_\infty\Bigl(\begin{psmallmatrix} x^*(T) \\ u^*(T) \end{psmallmatrix}\Bigr). \label{eq:bellman}
\end{align}

\begin{proof}
This follows as in~\citeA[Th.~9]{Ilch2018}.
\end{proof}
\end{proposition}

\section{MPC: asymptotic stability of the origin}\label{sec:stability}

To prove asymptotic stability of the origin w.r.t.\ the~MPC scheme from~\Cref{sec:mpc}, we employ the equivalent~ODE constructed in~\Cref{thm:odeequivalence}.

\begin{definition}
The set \(\mathbb{X}_\mathrm{f}\subseteq\mathbb{R}^n\) is called \emph{controlled forward invariant} w.r.t.\ the system \([E,A,B]\in\Sigma\) if, and only if,
\begin{equation*}
	\forall x^0\in\mathbb{X}_\mathrm{f}\ \exists(x,u)\in\mathfrak{B}_{[E,A,B]}\,\overset{\mathrm{aa}}{\forall} t\geq0: x(t)\in\mathbb{X}_\mathrm{f} \,\land\, E x^0=E x(0)
\end{equation*}
\end{definition}

The following theorem states that the optimal solution fulfils the condition in the preceding definition.

\begin{theorem}\label{thm:terminalregion}
Consider \([E,A,B]\in\Sigma\) with constraints~\eqref{eq:constraints}. 
Let the transformation matrix~\(\widehat{T}\) be defined as in~\eqref{eq:bergdoortrans}, %
and \(\widehat{B}\), \(\widehat{H}\), \(\widehat{R}\) 
 be defined by~\eqref{eq:abhat} and~\eqref{eq:shat}. Denote by~\(\widehat{P}\) the solution of the algebraic Riccati equation~\eqref{eq:are}. Define
\begin{equation*}
	\rho := \lambda_{\mathrm{min}}(\widehat{P}) \norm[\Bigg]{\begin{bmatrix} F & G \end{bmatrix} \widehat{T} \begin{bmatrix} I_{\hat{n}} \\ -\widehat{R}^{-1} (\widehat{B}^\top \widehat{P}+\widehat{H}) \end{bmatrix}}_\infty^{-2}>0.
\end{equation*}
Then the set
\begin{equation}
	\mathbb{X}_\mathrm{f} := \Biggl\{ [I_n,0] \widehat{T} \begin{bsmallmatrix} I_{\hat{n}} \\ -\widehat{R}^{-1} (\widehat{B}^\top \widehat{P}+\widehat{H}) \end{bsmallmatrix} \hat{x} \,\Biggm|\, \hat{x}\in\mathbb{R}^{\hat{n}} \,\land\, \hat{x}^\top \widehat{P} \hat{x} \leq \rho \Biggr\} \label{eq:terminalregion}
\end{equation}
is controlled forward invariant. Moreover, the optimal solution consisting of \(x^*\) and \(u^*\) satisfies the decrease condition~\eqref{eq:DecentCondition} almost everywhere.

\begin{proof}
Let~\(x^0\in\mathbb{X}_\mathrm{f}\) be arbitrary, and let~\((x,u)\in\mathfrak{B}_{[E,A,B]}\) be any solution with~\((E x)(0)=x^0\). Consider the solution of the~ODE
\begin{equation*}
	\dot{\hat{x}}(t) = [\widehat{A}-\widehat{B} \widehat{R}^{-1} (\widehat{B}^\top \widehat{P}+\widehat{H})] \hat{x}(t), \quad
	\hat{x}(0) = X^{-1} \begin{pmatrix} x^0 \\ u(0) \end{pmatrix}.
\end{equation*}
By the Bellman equation, it holds that
\begin{equation*}
	\forall t\geq0: \hat{x}(t)^\top \widehat{P} \hat{x}(t) \leq \hat{x}(0)^\top \widehat{P} \hat{x}(0) \leq \rho.
\end{equation*}
Note that~\(\hat{x}\in\mathcal{C}^\infty(\mathbb{R}_{\geq0},\mathbb{R}^r)\), therefore by~\Cref{thm:feedbackregularization} and~\Cref{thm:odeequivalence}, it holds that
\begin{align*}
	\begin{pmatrix} \tilde{x} \\ \tilde{u} \end{pmatrix} &= \widehat{T} \begin{bsmallmatrix} I_{\hat{n}} \\ -\widehat{R}^{-1} (\widehat{B}^\top \widehat{P}+\widehat{H}) \end{bsmallmatrix} x 
	\in \mathfrak{B}_{[E,A,B]}.
\end{align*}
For~\((\tilde{x},\tilde{u})\), it follows that~\((E \tilde{x})(0)=E x^0\) and for \(t\geq0\),
\begin{align*}
	\mathrel{\hphantom{=}} \norm{F \tilde{x}(t)+G \tilde{u}(t)}_\infty
	&= \norm[\Bigg]{\begin{bmatrix} F & G \end{bmatrix} \widehat{T} \begin{bmatrix} I_{\hat{n}} \\ -\widehat{R}^{-1} (\widehat{B}^\top \widehat{P}+\widehat{H}) \end{bmatrix} \hat{x}(t)}_\infty \\
	&\leq \norm[\Bigg]{\begin{bmatrix} F & G \end{bmatrix} \widehat{T} \begin{bmatrix} I_{\hat{n}} \\ -\widehat{R}^{-1} (\widehat{B}^\top \widehat{P}+\widehat{H}) \end{bmatrix}}_\infty \norm{\hat{x}(t)}_2 \\
	&\leq \norm[\Bigg]{\begin{bmatrix} F & G \end{bmatrix} \widehat{T}}_\infty \sqrt{\frac{\hat{x}(t) \widehat{P} \hat{x}(t)}{\lambda_{\mathrm{min}}(\widehat{P})}} \\
	&\leq 1.
\end{align*}
This shows that~\(\mathbb{X}_\mathrm{f}\) is controlled forward invariant.

Satisfaction of the decrease condition~\eqref{eq:DecentCondition} follows immediately for the optimal solution of the~DAE-OCP~\eqref{eq:ocp} (guaranteed to exist by~\Cref{thm:existenceoptimaltrajectory}). In light of the Bellman equation~\eqref{eq:bellman}, the terminal cost is simply given by the optimal cost~\(V_\infty\).
\end{proof}
\end{theorem}

\section{Example}\label{sec:example}
\begin{figure}
	\includegraphics[width=\columnwidth]{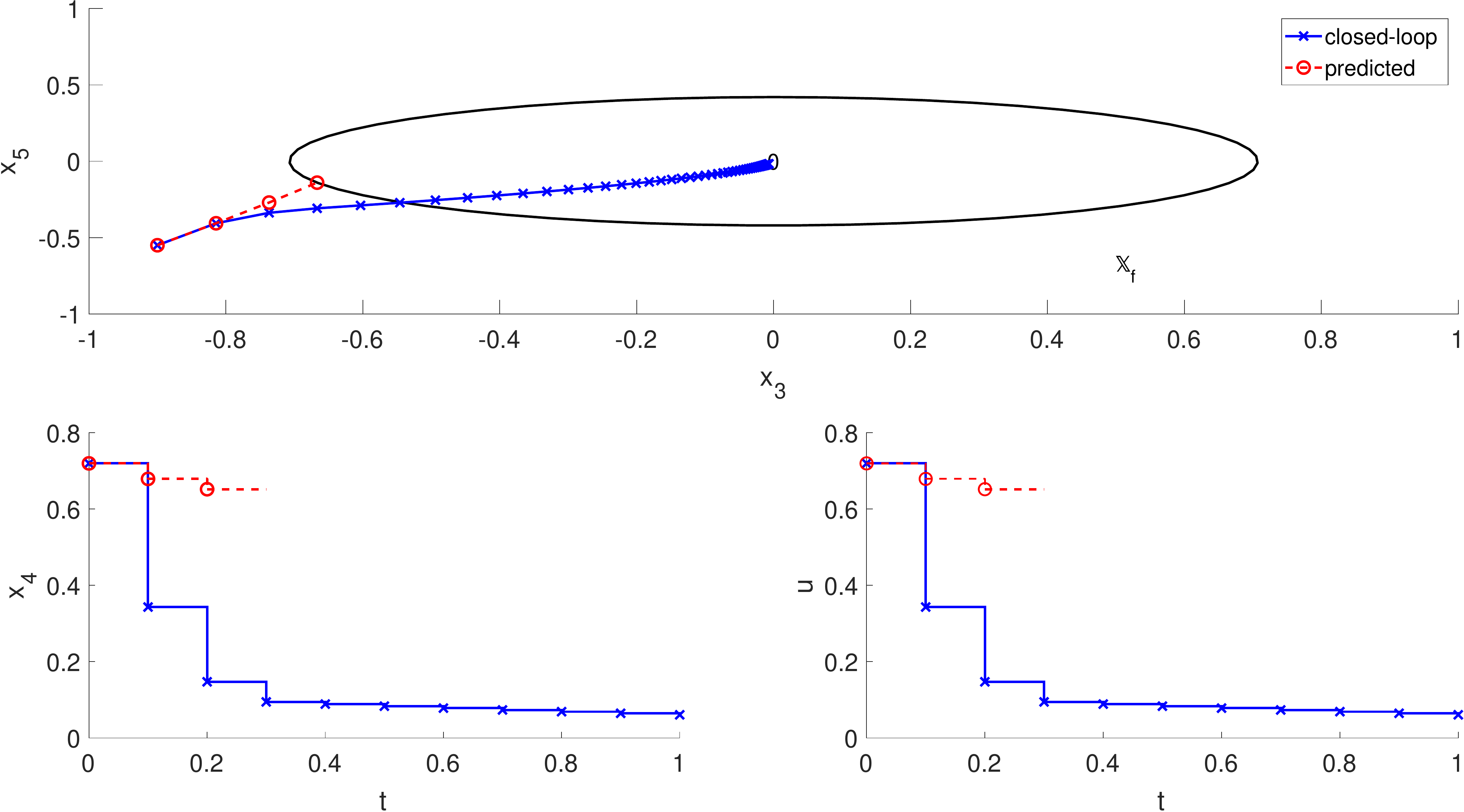}
	\caption{Closed-loop performance of the MPC scheme and predicted solution of the optimal control problem with added terminal constraints and costs from \ref{item:step2} of the MPC algorithm at time \(t=k \delta=0\)}\label{fig:example}
\end{figure}

\noindent
Minimize the cost functional
\begin{equation}\nonumber
	\int_0^T \norm{x(t)}^2+\norm{u(t)}^2 \mathrm{d} t
\end{equation}
subject to the singular~DAE
\begin{equation}\nonumber
	\frac{\mathrm{d}}{\mathrm{d} t} \left[\begin{array}{r|r|r|rr} 0 & 0 & 0 & 0 & 0 \\ 1 & 0 & 0 & 0 & 0 \\ \hline 0 & 0 & 0 & 0 & 0 \\ \hline 0 & 0 & 1 & 0 & 0 \\ \hline 0 & 0 & 0 & 0 & 1 \end{array}\right] x(t) = \left[\begin{array}{r|r|r|rr} 1 & 0 & 0 & 0 & 0 \\ 0 & 0 & 0 & 0 & 0 \\ \hline 0 & 1 & 0 & 0 & 0 \\ \hline 0 & 0 & -1 & 0 & 0 \\ \hline 0 & 0 & 0 & 1 & 0 \end{array}\right] x(t) + \begin{pmatrix} 0 \\ 0 \\ 0 \\ 0 \\ 1 \end{pmatrix} u(t)
\end{equation}
and the initial condition
\begin{equation}\nonumber
	(E x)(0) = E x^0.
\end{equation}
The~ODE obtained from regularization is given by
\begin{equation}\label{eq:exampleequivalentode}
	\dot{z}_1(t) = \begin{bmatrix} -1 & 0 \\ 0 & 0 \end{bmatrix} z_1(t) + \begin{bmatrix}  0 & 0 \\ 1 & 1 \end{bmatrix} v(t), \quad z_2(t)=0, 
\end{equation}
where \(z_{1,1}=x_3\), \(z_{1,2}=x_5\), \(z_{2,1}=x_1\), \(z_{2,2}=-x_2\), \(v_1=x_4\), \(v_2=u\), and the equivalent~OCP is
\begin{align*}
	&\mathrm{Minimize} \int_0^T \norm{z_1(t)}^2+\norm{v(t)}^2 \mathrm{d} t \\
	&\text{s.t.~\eqref{eq:exampleequivalentode} with } z_1(0) = \begin{pmatrix} x^0_3 \\ x^0_5 \end{pmatrix}.
\end{align*}
Obviously~\Cref{assumption} is satisfied. 
The solution of the algebraic Riccati equation~\eqref{eq:are} is \(\diag(\tfrac{1}{2} ,\,\sqrt{2})\). For the constraints
\begin{equation*}
	-1\leq x_i(t) \leq 1, i\in\{1,\dots,5\}, \quad -1\leq u(t)\leq 1,
\end{equation*}
written in matrix form as
\begin{equation*}
	\begin{bmatrix*}[r] I_5 & 0 \\ -I_5 & 0 \\ 0 & 1 \\ 0 & -1 \end{bmatrix*} \begin{pmatrix} x(t) \\ u(t) \end{pmatrix} \leq \mathds{1},
\end{equation*}
we obtain
\begin{equation*}
\begin{lgathered}
	\rho = \frac{1}{2} \norm*{\begin{bmatrix} 0 & 0 \\ 0 & 0 \\ 1 & 0 \\ 0 & \sqrt{2} \\ 0 & 1 \end{bmatrix}}_\infty^{-2} = \frac{1}{4},
	\quad
	\mathbb{X}_\mathrm{f} = \left\{ \begin{bmatrix} 0 & 0 \\ 0 & 0 \\ 1 & 0 \\ 0 & \sqrt{2} \\ 0 & 1 \end{bmatrix} \hat{x}  \,\middle|\, \hat{x}\in\mathbb{R}^2 \,\land\, \frac{1}{2} \hat{x}_1^2+\sqrt{2} \hat{x}_2^2 \leq \frac{1}{4}\right\}, \\
	V_\mathrm{f}(\hat{x}) = \frac{1}{2} \hat{x}_3^2+\sqrt{2} \hat{x}_5^2, \quad \hat{x}\in\mathbb{X}_\mathrm{f}.
\end{lgathered}
\end{equation*}
The constructed terminal region~\(\mathbb{X}_\mathrm{f}\) and the performance for a~MPC scheme with step size~\(\delta=0.1\) and prediction horizon~\(T=3\delta\) is depicted in~\Cref{fig:example}. The states~\(x_1\) and~\(x_2\) are omitted as it follows from the~DAE that~\(x_1=0\), \(x_2\overset{\text{ae}}{=}0\). Starting with an initial value of \((E x)(0)=(0,0,0-0.9,-0.55)\), the closed-loop solution for~\(x\) converges to the origin. In addition, the figure also depicts the solution of the optimal control problem with added terminal constraints and costs from \ref{item:step2} of the MPC algorithm at time \(t=k \delta=0\). It can be seen that the state reaches the boundary of the terminal region within the prediction horizon \(T=3\delta\), since that is mandated by the terminal constraint. On the other hand, the actual closed-loop MPC solution does not reach the terminal region within this time due to the receding-horizon nature of the MPC scheme. Once the interior of the terminal region is reached by the closed-loop solution, it is never left again. This follows since the MPC closed-loop solution coincides with the (unconstrained) infinite-horizon optimal solution by construction once  the terminal region is reached.

\section{Conclusions and open problem}
In this paper we describe a way to obtain a~MPC scheme for a~DAE with state and input constraints that guarantees stability of the closed loop w.r.t.\ the origin. This is achieved by regularizing the~DAE to obtain an~ODE optimal control problem for which a terminal region and costs can be constructed. These terminal ingredients can then be expressed in terms of the nominal~DAE by a state transformation.

In the future, we want to investigate whether it is possible to achieve similar results without having to resort to a transformation to an equivalent~ODE. This would allow to express~\Cref{assumption} directly in terms of the~DAE: for example, it can be easily seen that the stabilizability of the equivalent~ODE 
in~\Cref{assumption} is equivalent to the behavioural stabilizability of the nominal~DAE~\eqref{eq:dae}. We would like to obtain similar results for the rest of the assumptions where the situation is much less obvious.

One way to work directly with the~DAE is to adapt the approach 
by~\citeA{ReisVoig18} to model predictive control: their results allow to characterize the optimal value and optimal solution using so-called Lur'e equations for the~DAE. In order to use these findings for~MPC, it is necessary to characterize the positive definiteness of the optimal value in terms of the~DAE-OCP. Using these results, a construction similar to~\Cref{thm:terminalregion} yields a terminal region that, together with the optimal value as terminal costs, guarantees asymptotic stability of the~MPC scheme w.r.t.\  the origin.

\bigskip

\textbf{Acknowledgements}: The authors are indebted to the German Research Foundation~(DFG) (grants~IL25/10-1, WO2056/2-1, RE2917/4-1, and WO2056/6-1) and the Studienstiftung
des Deutschen Volkes for their support. Furthermore we thank our colleague Thomas Berger (Paderborn) for constructive suggestions and discussions.


\begin{thebibliography}{}

\bibitem [\protect \citeauthoryear {%
Bankmann%
}{%
Bankmann%
}{%
{\protect \APACyear {2016}}%
}]{%
bankmann2016linear}
\APACinsertmetastar {%
bankmann2016linear}%
\begin{APACrefauthors}%
Bankmann, D.%
\end{APACrefauthors}%
\unskip\
\newblock
\APACrefYear{2016}.
\unskip\
\newblock
\APACrefbtitle {On linear-quadratic control theory of implicit difference
  equations} {On linear-quadratic control theory of implicit difference
  equations}\ \APACtypeAddressSchool {\BUMTh}{}{}.
\unskip\
\newblock
\APACaddressSchool {Berlin, Germany}{Fakult\"at II - Mathematik und
  Naturwissenschaften, Technische Universit{\"a}t Berlin}.
\PrintBackRefs{\CurrentBib}

\bibitem [\protect \citeauthoryear {%
Benner%
, Losse%
, Mehrmann%
\BCBL {}\ \BBA {} Voigt%
}{%
Benner%
\ \protect \BOthers {.}}{%
{\protect \APACyear {2015}}%
}]{%
BennLoss15}
\APACinsertmetastar {%
BennLoss15}%
\begin{APACrefauthors}%
Benner, P.%
, Losse, P.%
, Mehrmann, V.%
\BCBL {}\ \BBA {} Voigt, M.%
\end{APACrefauthors}%
\unskip\
\newblock
\APACrefYearMonthDay{2015}{}{}.
\newblock
{\BBOQ}\APACrefatitle {Numerical Linear Algebra Methods for Linear
  Differential-Algebraic Equations} {Numerical linear algebra methods for
  linear differential-algebraic equations}.{\BBCQ}
\newblock
\BIn{} A.~Ilchmann\ \BBA {} T.~Reis\ (\BEDS), \APACrefbtitle {Surveys in
  Differential-Algebraic Equations {III}} {Surveys in differential-algebraic
  equations {III}}\ (\BPG~117-175).
\newblock
\APACaddressPublisher{Berlin-Heidelberg}{Springer-Verlag}.
\PrintBackRefs{\CurrentBib}

\bibitem [\protect \citeauthoryear {%
Berger%
\ \BBA {} Reis%
}{%
Berger%
\ \BBA {} Reis%
}{%
{\protect \APACyear {2013}}%
}]{%
BergReis13a}
\APACinsertmetastar {%
BergReis13a}%
\begin{APACrefauthors}%
Berger, T.%
\BCBT {}\ \BBA {} Reis, T.%
\end{APACrefauthors}%
\unskip\
\newblock
\APACrefYearMonthDay{2013}{}{}.
\newblock
{\BBOQ}\APACrefatitle {Controllability of linear differential-algebraic
  systems~--~a survey} {Controllability of linear differential-algebraic
  systems~--~a survey}.{\BBCQ}
\newblock
\BIn{} A.~Ilchmann\ \BBA {} T.~Reis\ (\BEDS), \APACrefbtitle {Surveys in
  Differential-Algebraic Equations {I}} {Surveys in differential-algebraic
  equations {I}}\ (\BPGS\ 1--61).
\newblock
\APACaddressPublisher{Berlin-Heidelberg}{Springer-Verlag}.
\PrintBackRefs{\CurrentBib}

\bibitem [\protect \citeauthoryear {%
Berger%
\ \BBA {} Trenn%
}{%
Berger%
\ \BBA {} Trenn%
}{%
{\protect \APACyear {2012}}%
}]{%
BergTren12}
\APACinsertmetastar {%
BergTren12}%
\begin{APACrefauthors}%
Berger, T.%
\BCBT {}\ \BBA {} Trenn, S.%
\end{APACrefauthors}%
\unskip\
\newblock
\APACrefYearMonthDay{2012}{}{}.
\newblock
{\BBOQ}\APACrefatitle {The quasi-{K}ronecker form for matrix pencils} {The
  quasi-{K}ronecker form for matrix pencils}.{\BBCQ}
\newblock
\APACjournalVolNumPages{{SIAM} J. Matrix Anal. \& Appl.}{33}{2}{336--368}.
\PrintBackRefs{\CurrentBib}

\bibitem [\protect \citeauthoryear {%
Berger%
\ \BBA {} Van~Dooren%
}{%
Berger%
\ \BBA {} Van~Dooren%
}{%
{\protect \APACyear {2015}}%
}]{%
BergDoor15}
\APACinsertmetastar {%
BergDoor15}%
\begin{APACrefauthors}%
Berger, T.%
\BCBT {}\ \BBA {} Van~Dooren, P.%
\end{APACrefauthors}%
\unskip\
\newblock
\APACrefYearMonthDay{2015}{}{}.
\newblock
{\BBOQ}\APACrefatitle {Computing the regularization of a linear
  differential-algebraic system} {Computing the regularization of a linear
  differential-algebraic system}.{\BBCQ}
\newblock
\APACjournalVolNumPages{Syst. Control Lett.}{86}{}{48--53}.
\PrintBackRefs{\CurrentBib}

\bibitem [\protect \citeauthoryear {%
Bock%
, Diehl%
, Kostina%
\BCBL {}\ \BBA {} Schl{\"o}der%
}{%
Bock%
\ \protect \BOthers {.}}{%
{\protect \APACyear {2007}}%
}]{%
bock2007constrained}
\APACinsertmetastar {%
bock2007constrained}%
\begin{APACrefauthors}%
Bock, H\BPBI G.%
, Diehl, M.%
, Kostina, E.%
\BCBL {}\ \BBA {} Schl{\"o}der, J\BPBI P.%
\end{APACrefauthors}%
\unskip\
\newblock
\APACrefYearMonthDay{2007}{}{}.
\newblock
{\BBOQ}\APACrefatitle {Constrained optimal feedback control of systems governed
  by large differential algebraic equations} {Constrained optimal feedback
  control of systems governed by large differential algebraic
  equations}.{\BBCQ}
\newblock
\BIn{} L\BPBI T.~Biegler, O.~Ghattas, M.~Heinkenschloss, D.~Keyes\BCBL {}\ \BBA
  {} B.~van Bloemen~Waanders\ (\BEDS), \APACrefbtitle {Real-Time
  {PDE}-Constrained Optimization} {Real-time {PDE}-constrained optimization}\
  (\BPGS\ 3--24).
\newblock
\APACaddressPublisher{Philadelphia}{SIAM}.
\PrintBackRefs{\CurrentBib}

\bibitem [\protect \citeauthoryear {%
Bunse-Gerstner%
, Mehrmann%
\BCBL {}\ \BBA {} Nichols%
}{%
Bunse-Gerstner%
\ \protect \BOthers {.}}{%
{\protect \APACyear {1992}}%
}]{%
BunsMehr92}
\APACinsertmetastar {%
BunsMehr92}%
\begin{APACrefauthors}%
Bunse-Gerstner, A.%
, Mehrmann, V.%
\BCBL {}\ \BBA {} Nichols, N\BPBI K.%
\end{APACrefauthors}%
\unskip\
\newblock
\APACrefYearMonthDay{1992}{}{}.
\newblock
{\BBOQ}\APACrefatitle {Regularization of descriptor systems by derivative and
  proportional state feedback} {Regularization of descriptor systems by
  derivative and proportional state feedback}.{\BBCQ}
\newblock
\APACjournalVolNumPages{{SIAM} J. Matrix Anal. \& Appl.}{13}{1}{46--67}.
\PrintBackRefs{\CurrentBib}

\bibitem [\protect \citeauthoryear {%
Bunse-Gerstner%
, Mehrmann%
\BCBL {}\ \BBA {} Nichols%
}{%
Bunse-Gerstner%
\ \protect \BOthers {.}}{%
{\protect \APACyear {1994}}%
}]{%
BunsMehr94}
\APACinsertmetastar {%
BunsMehr94}%
\begin{APACrefauthors}%
Bunse-Gerstner, A.%
, Mehrmann, V.%
\BCBL {}\ \BBA {} Nichols, N\BPBI K.%
\end{APACrefauthors}%
\unskip\
\newblock
\APACrefYearMonthDay{1994}{}{}.
\newblock
{\BBOQ}\APACrefatitle {Regularization of Descriptor Systems by Output Feedback}
  {Regularization of descriptor systems by output feedback}.{\BBCQ}
\newblock
\APACjournalVolNumPages{{IEEE} Trans. Autom. Control}{39}{8}{1742--1748}.
\PrintBackRefs{\CurrentBib}

\bibitem [\protect \citeauthoryear {%
S.~Campbell%
, Ilchmann%
, Mehrmann%
\BCBL {}\ \BBA {} Reis%
}{%
S.~Campbell%
\ \protect \BOthers {.}}{%
{\protect \APACyear {2019}}%
}]{%
campbellapplications}
\APACinsertmetastar {%
campbellapplications}%
\begin{APACrefauthors}%
Campbell, S.%
, Ilchmann, A.%
, Mehrmann, V.%
\BCBL {}\ \BBA {} Reis, T.%
\end{APACrefauthors}%
\ (\BEDS).
\unskip\
\newblock
\APACrefYear{2019}.
\newblock
\APACrefbtitle {Applications of Differential-Algebraic Equations: Examples and
  Benchmarks} {Applications of differential-algebraic equations: Examples and
  benchmarks}.
\newblock
\APACaddressPublisher{Cham}{Springer}.
\PrintBackRefs{\CurrentBib}

\bibitem [\protect \citeauthoryear {%
S\BPBI L.~Campbell%
, Kunkel%
\BCBL {}\ \BBA {} Mehrmann%
}{%
S\BPBI L.~Campbell%
\ \protect \BOthers {.}}{%
{\protect \APACyear {2012}}%
}]{%
CampKunk12}
\APACinsertmetastar {%
CampKunk12}%
\begin{APACrefauthors}%
Campbell, S\BPBI L.%
, Kunkel, P.%
\BCBL {}\ \BBA {} Mehrmann, V.%
\end{APACrefauthors}%
\unskip\
\newblock
\APACrefYearMonthDay{2012}{}{}.
\newblock
{\BBOQ}\APACrefatitle {Regularization of Linear and Nonlinear Descriptor
  Systems} {Regularization of linear and nonlinear descriptor systems}.{\BBCQ}
\newblock
\BIn{} L\BPBI T.~Biegler, S\BPBI L.~Campbell\BCBL {}\ \BBA {} V.~Mehrmann\
  (\BEDS), \APACrefbtitle {Control and Optimization with Differential-Algebraic
  Constraints} {Control and optimization with differential-algebraic
  constraints}\ (\BPGS\ 17--36).
\newblock
\APACaddressPublisher{Philadelphia}{SIAM}.
\PrintBackRefs{\CurrentBib}

\bibitem [\protect \citeauthoryear {%
Cobb%
}{%
Cobb%
}{%
{\protect \APACyear {1983}}%
}]{%
Cobb83}
\APACinsertmetastar {%
Cobb83}%
\begin{APACrefauthors}%
Cobb, J\BPBI D.%
\end{APACrefauthors}%
\unskip\
\newblock
\APACrefYearMonthDay{1983}{}{}.
\newblock
{\BBOQ}\APACrefatitle {Descriptor variable systems and optimal state
  regulation} {Descriptor variable systems and optimal state
  regulation}.{\BBCQ}
\newblock
\APACjournalVolNumPages{{IEEE} Trans. Autom. Control}{28}{}{601--611}.
\PrintBackRefs{\CurrentBib}

\bibitem [\protect \citeauthoryear {%
Coron%
, Gr\"{u}ne%
\BCBL {}\ \BBA {} Worthmann%
}{%
Coron%
\ \protect \BOthers {.}}{%
{\protect \APACyear {2020}}%
}]{%
CoroGrun20}
\APACinsertmetastar {%
CoroGrun20}%
\begin{APACrefauthors}%
Coron, J\BHBI M.%
, Gr\"{u}ne, L.%
\BCBL {}\ \BBA {} Worthmann, K.%
\end{APACrefauthors}%
\unskip\
\newblock
\APACrefYearMonthDay{2020}{}{}.
\newblock
{\BBOQ}\APACrefatitle {Model Predictive Control, Cost Controllability, and
  Homogeneity} {Model predictive control, cost controllability, and
  homogeneity}.{\BBCQ}
\newblock
\APACjournalVolNumPages{SIAM J. Control Optim.}{58}{5}{2979--2996}.
\PrintBackRefs{\CurrentBib}

\bibitem [\protect \citeauthoryear {%
Diehl%
\ \protect \BOthers {.}}{%
Diehl%
\ \protect \BOthers {.}}{%
{\protect \APACyear {2002}}%
}]{%
DiehBock02}
\APACinsertmetastar {%
DiehBock02}%
\begin{APACrefauthors}%
Diehl, M.%
, Bock, H\BPBI G.%
, Schl{\"o}der, J\BPBI P.%
, Findeisen, R.%
, Nagy, Z.%
\BCBL {}\ \BBA {} Allg{\"o}wer, F.%
\end{APACrefauthors}%
\unskip\
\newblock
\APACrefYearMonthDay{2002}{}{}.
\newblock
{\BBOQ}\APACrefatitle {Real-time optimization and nonlinear model predictive
  control of processes governed by differential-algebraic equations} {Real-time
  optimization and nonlinear model predictive control of processes governed by
  differential-algebraic equations}.{\BBCQ}
\newblock
\APACjournalVolNumPages{Journal of {P}rocess {C}ontrol}{12}{}{577--585}.
\PrintBackRefs{\CurrentBib}

\bibitem [\protect \citeauthoryear {%
Gerdts%
}{%
Gerdts%
}{%
{\protect \APACyear {2011}}%
}]{%
gerdts2011optimal}
\APACinsertmetastar {%
gerdts2011optimal}%
\begin{APACrefauthors}%
Gerdts, M.%
\end{APACrefauthors}%
\unskip\
\newblock
\APACrefYear{2011}.
\newblock
\APACrefbtitle {Optimal control of {ODE}s and {DAE}s} {Optimal control of
  {ODE}s and {DAE}s}.
\newblock
\APACaddressPublisher{Berlin}{Walter de Gruyter}.
\PrintBackRefs{\CurrentBib}

\bibitem [\protect \citeauthoryear {%
Ilchmann%
, Leben%
, Witschel%
\BCBL {}\ \BBA {} Worthmann%
}{%
Ilchmann%
\ \protect \BOthers {.}}{%
{\protect \APACyear {2019}}%
}]{%
ilchmann2019optimal}
\APACinsertmetastar {%
ilchmann2019optimal}%
\begin{APACrefauthors}%
Ilchmann, A.%
, Leben, L.%
, Witschel, J.%
\BCBL {}\ \BBA {} Worthmann, K.%
\end{APACrefauthors}%
\unskip\
\newblock
\APACrefYearMonthDay{2019}{}{}.
\newblock
{\BBOQ}\APACrefatitle {Optimal control of differential-algebraic equations from
  an ordinary differential equation perspective} {Optimal control of
  differential-algebraic equations from an ordinary differential equation
  perspective}.{\BBCQ}
\newblock
\APACjournalVolNumPages{Optimal Control Applications and
  Methods}{40}{2}{351--366}.
\PrintBackRefs{\CurrentBib}

\bibitem [\protect \citeauthoryear {%
Ilchmann%
, Witschel%
\BCBL {}\ \BBA {} Worthmann%
}{%
Ilchmann%
\ \protect \BOthers {.}}{%
{\protect \APACyear {2018}}%
}]{%
Ilch2018}
\APACinsertmetastar {%
Ilch2018}%
\begin{APACrefauthors}%
Ilchmann, A.%
, Witschel, J.%
\BCBL {}\ \BBA {} Worthmann, K.%
\end{APACrefauthors}%
\unskip\
\newblock
\APACrefYearMonthDay{2018}{}{}.
\newblock
{\BBOQ}\APACrefatitle {Model Predictive Control for Linear
  Differential-Algebraic Equations} {Model predictive control for linear
  differential-algebraic equations}.{\BBCQ}
\newblock
\APACjournalVolNumPages{IFAC PapersOnLine}{51}{}{98-103}.
\newblock
\APACrefnote{6th IFAC Conference on Nonlinear Model Predictive Control,
  Madison, Wisconsin, USA}
\PrintBackRefs{\CurrentBib}

\bibitem [\protect \citeauthoryear {%
Kouvaritakis%
\ \BBA {} Cannon%
}{%
Kouvaritakis%
\ \BBA {} Cannon%
}{%
{\protect \APACyear {2016}}%
}]{%
kouvaritakis2016model}
\APACinsertmetastar {%
kouvaritakis2016model}%
\begin{APACrefauthors}%
Kouvaritakis, B.%
\BCBT {}\ \BBA {} Cannon, M.%
\end{APACrefauthors}%
\unskip\
\newblock
\APACrefYear{2016}.
\newblock
\APACrefbtitle {Model Predictive Control} {Model predictive control}.
\newblock
\APACaddressPublisher{Cham}{Springer}.
\PrintBackRefs{\CurrentBib}

\bibitem [\protect \citeauthoryear {%
Kunkel%
\ \BBA {} Mehrmann%
}{%
Kunkel%
\ \BBA {} Mehrmann%
}{%
{\protect \APACyear {2008}}%
}]{%
KunkMehr08}
\APACinsertmetastar {%
KunkMehr08}%
\begin{APACrefauthors}%
Kunkel, P.%
\BCBT {}\ \BBA {} Mehrmann, V.%
\end{APACrefauthors}%
\unskip\
\newblock
\APACrefYearMonthDay{2008}{}{}.
\newblock
{\BBOQ}\APACrefatitle {Optimal control for unstructured nonlinear
  differential-algebraic equations of arbitrary index} {Optimal control for
  unstructured nonlinear differential-algebraic equations of arbitrary
  index}.{\BBCQ}
\newblock
\APACjournalVolNumPages{Math. Control Signals Syst.}{20}{}{227--269}.
\PrintBackRefs{\CurrentBib}

\bibitem [\protect \citeauthoryear {%
Lamour%
, M{\"a}rz%
\BCBL {}\ \BBA {} Tischendorf%
}{%
Lamour%
\ \protect \BOthers {.}}{%
{\protect \APACyear {2013}}%
}]{%
LamoMarz13}
\APACinsertmetastar {%
LamoMarz13}%
\begin{APACrefauthors}%
Lamour, R.%
, M{\"a}rz, R.%
\BCBL {}\ \BBA {} Tischendorf, C.%
\end{APACrefauthors}%
\unskip\
\newblock
\APACrefYear{2013}.
\newblock
\APACrefbtitle {Differential-Algebraic Equations: A Projector Based Analysis}
  {Differential-algebraic equations: A projector based analysis}.
\newblock
\APACaddressPublisher{Heidelberg-Berlin}{Springer-Verlag}.
\PrintBackRefs{\CurrentBib}

\bibitem [\protect \citeauthoryear {%
Lancaster%
\ \BBA {} Rodman%
}{%
Lancaster%
\ \BBA {} Rodman%
}{%
{\protect \APACyear {1995}}%
}]{%
LancRodm95}
\APACinsertmetastar {%
LancRodm95}%
\begin{APACrefauthors}%
Lancaster, P.%
\BCBT {}\ \BBA {} Rodman, L.%
\end{APACrefauthors}%
\unskip\
\newblock
\APACrefYear{1995}.
\newblock
\APACrefbtitle {Algebraic {R}iccati Equations} {Algebraic {R}iccati equations}.
\newblock
\APACaddressPublisher{Oxford}{Clarendon Press}.
\PrintBackRefs{\CurrentBib}

\bibitem [\protect \citeauthoryear {%
Mayne%
, Rawlings%
, Rao%
\BCBL {}\ \BBA {} Scokaert%
}{%
Mayne%
\ \protect \BOthers {.}}{%
{\protect \APACyear {2000}}%
}]{%
Mayne2000}
\APACinsertmetastar {%
Mayne2000}%
\begin{APACrefauthors}%
Mayne, D\BPBI Q.%
, Rawlings, J\BPBI B.%
, Rao, C\BPBI V.%
\BCBL {}\ \BBA {} Scokaert, P\BPBI O.%
\end{APACrefauthors}%
\unskip\
\newblock
\APACrefYearMonthDay{2000}{}{}.
\newblock
{\BBOQ}\APACrefatitle {Constrained model predictive control: Stability and
  optimality} {Constrained model predictive control: Stability and
  optimality}.{\BBCQ}
\newblock
\APACjournalVolNumPages{Automatica}{36}{6}{789-814}.
\PrintBackRefs{\CurrentBib}

\bibitem [\protect \citeauthoryear {%
Rawlings%
, Mayne%
\BCBL {}\ \BBA {} Diehl%
}{%
Rawlings%
\ \protect \BOthers {.}}{%
{\protect \APACyear {2017}}%
}]{%
RawlingsMayneDiehl2017}
\APACinsertmetastar {%
RawlingsMayneDiehl2017}%
\begin{APACrefauthors}%
Rawlings, J\BPBI B.%
, Mayne, D\BPBI Q.%
\BCBL {}\ \BBA {} Diehl, M\BPBI M.%
\end{APACrefauthors}%
\unskip\
\newblock
\APACrefYear{2017}.
\newblock
\APACrefbtitle {Model Predictive Control: Theory, Computation, and Design}
  {Model predictive control: Theory, computation, and design}\
  (\PrintOrdinal{2nd}\ \BEd).
\newblock
\APACaddressPublisher{Madison, Wisconsin}{Nob Hill Publishing}.
\PrintBackRefs{\CurrentBib}

\bibitem [\protect \citeauthoryear {%
Reis%
\ \BBA {} Voigt%
}{%
Reis%
\ \BBA {} Voigt%
}{%
{\protect \APACyear {2019}}%
}]{%
ReisVoig18}
\APACinsertmetastar {%
ReisVoig18}%
\begin{APACrefauthors}%
Reis, T.%
\BCBT {}\ \BBA {} Voigt, M.%
\end{APACrefauthors}%
\unskip\
\newblock
\APACrefYearMonthDay{2019}{}{}.
\newblock
{\BBOQ}\APACrefatitle {Linear-quadratic optimal control of
  differential-algebraic systems: {T}he infinite time horizon problem with zero
  terminal state} {Linear-quadratic optimal control of differential-algebraic
  systems: {T}he infinite time horizon problem with zero terminal
  state}.{\BBCQ}
\newblock
\APACjournalVolNumPages{{SIAM} J. Control Optim.}{57}{3}{1567--1596}.
\PrintBackRefs{\CurrentBib}

\bibitem [\protect \citeauthoryear {%
Riaza%
}{%
Riaza%
}{%
{\protect \APACyear {2008}}%
}]{%
Riaz08}
\APACinsertmetastar {%
Riaz08}%
\begin{APACrefauthors}%
Riaza, R.%
\end{APACrefauthors}%
\unskip\
\newblock
\APACrefYear{2008}.
\newblock
\APACrefbtitle {Differential-Algebraic Systems: {A}nalytical Aspects and
  Circuit Applications} {Differential-algebraic systems: {A}nalytical aspects
  and circuit applications}.
\newblock
\APACaddressPublisher{Basel}{World Scientific Publishing}.
\PrintBackRefs{\CurrentBib}

\bibitem [\protect \citeauthoryear {%
Sj{\"o}berg%
, Findeisen%
\BCBL {}\ \BBA {} Allg{\"o}wer%
}{%
Sj{\"o}berg%
\ \protect \BOthers {.}}{%
{\protect \APACyear {2007}}%
}]{%
SjobFind07}
\APACinsertmetastar {%
SjobFind07}%
\begin{APACrefauthors}%
Sj{\"o}berg, J.%
, Findeisen, R.%
\BCBL {}\ \BBA {} Allg{\"o}wer, F.%
\end{APACrefauthors}%
\unskip\
\newblock
\APACrefYearMonthDay{2007}{}{}.
\newblock
{\BBOQ}\APACrefatitle {Model Predictive Control of Continuous Time Nonlinear
  Differential Algebraic Systems} {Model predictive control of continuous time
  nonlinear differential algebraic systems}.{\BBCQ}
\newblock
\APACjournalVolNumPages{IFAC Proceedings Volumes}{40}{12}{48--53}.
\PrintBackRefs{\CurrentBib}

\bibitem [\protect \citeauthoryear {%
Yonchev%
, Findeisen%
, Ebenbauer%
\BCBL {}\ \BBA {} Allg{\"o}wer%
}{%
Yonchev%
\ \protect \BOthers {.}}{%
{\protect \APACyear {2004}}%
}]{%
YoncFind04}
\APACinsertmetastar {%
YoncFind04}%
\begin{APACrefauthors}%
Yonchev, A.%
, Findeisen, R.%
, Ebenbauer, C.%
\BCBL {}\ \BBA {} Allg{\"o}wer, F.%
\end{APACrefauthors}%
\unskip\
\newblock
\APACrefYearMonthDay{2004}{}{}.
\newblock
{\BBOQ}\APACrefatitle {Model predictive control of linear continuous time
  singular systems subject to input constraints} {Model predictive control of
  linear continuous time singular systems subject to input constraints}.{\BBCQ}
\newblock
\APACjournalVolNumPages{Proc. 43rd IEEE Conf. Decision Control
  (CDC)}{2}{}{2047--2052}.
\PrintBackRefs{\CurrentBib}

\end{thebibliography}
\end{document}